\documentclass{amsart}\address{Université Marie et Louis Pasteur, CNRS, LmB (UMR 6623), 25000 Besançon, France.}\author{Rubén Muñoz-\relax-Bertrand}\binoppenalty\maxdimen\relpenalty\maxdimen\title{Isocrystals and de Rham--Witt connections}\usepackage{amssymb,mathtools,mleftright,smfthm,stmaryrd,tikz-cd}\begin{document}\begin{abstract}We
introduce the notion of integrable connections for a sheaf of differential graded algebras on a topological space.
We then describe them in the finite locally projective setting, when the sheaf is either the de Rham complex of a formal or a weakly formal scheme, or for the convergent or the overconvergent de Rham--Witt complex on a smooth scheme over a perfect field of positive characteristic.
This enables us to give a new description of convergent and overconvergent isocrystals with a Frobenius
structure.\end{abstract}\maketitle\section*{Introduction}The
search for a $p$-adic Weil cohomology theory for schemes over a field of positive characteristic has a long and rich history.
The construction of Monsky--Washnitzer cohomology \cite{formalcohomologyi}, which was the earliest attempt to find such a theory, was inspired from the fact that on a smooth manifold, de Rham cohomology computes singular cohomology.
This allowed them to find a formula à la Lefschetz computing the zeta function of an affine scheme, under some assumptions.

In positive characteristic, the de Rham complex can have infinite dimensional cohomology groups.
In order to get a computable trace formula, one wants to work with finite dimensional vector spaces.
The idea of Monsky and Washnitzer was to lift the de Rham complex to another de Rham complex in characteristic $0$ having finite dimensional cohomology.
Moreover, they also lift the Frobenius endomorphism to get an action on these cohomology groups.

This action is paramount to the theory.
However, it is not possible in general to lift the absolute Frobenius of a scheme in positive characteristic to characteristic $0$; one can only do this locally.

In order to somehow glue these Frobenius lift globally, the first approach was crystalline cohomology \cite{cohomologiecristallinedesschemasdecaracteristique}.
The idea here is to consider a Grothendieck topology on the scheme, giving raise to the crystalline site.
The category of crystals are the coefficients for this cohomology, whose good properties hold for proper schemes.

A more general theory is rigid cohomology \cite{cohomologierigideetcohomologierigideasupportsproprespremierepartie}.
It retrieves both Monsky--Washnitzer and crystalline cohomology.
Berthelot's original construction did not use a site, but rather relied on constructions using rigid geometry.
It is now known that there actually is an overconvergent site computing rigid cohomology \cite{theoverconvergentsite}.

There have been other strategies to understand the global Frobenius action.
One can cite, for instance, the motivic approach which enables one to remove all the choices involved \cite{themonskywashnitzerandtheoverconvergentrealizations}.

In this article, we are interested in another viewpoint.
Crystalline and rigid cohomology can be computed without introducing any Grothendieck topology.
By using $p$-adic analytic methods, one can instead consider sheaves of differential graded algebras on the underlying topological space of the scheme.

This approach has been started with the introduction of the de Rham--Witt complex \cite{complexedederhamwittetcohomologiecristalline}.
This complex yields a quasi-isomorphism with crystalline cohomology.
Moreover, there exists a canonical and global lift of the Frobenius on the de Rham--Witt complex.
Also, by restricting to locally finite free crystals, one can see that they are also coefficients for this cohomology theory \cite{complexedederhamwittacoefficientsdansuncristal}\cite{crystalsandderhamwittconnections}.

There have been two different proposals for the computation of rigid cohomology.
The first one, the celebrated theory of arithmetic $\mathcal D$-modules, yields a category stable under the $6$-functor formalism.
Nevertheless, there are no global Frobenius lifts in that setting, so one has to resort to glueing to get the action.
Another viewpoint is the overconvergent de Rham--Witt complex \cite{overconvergentderhamwittcohomology}, which is known to compute rigid cohomology, but which had no theory of coefficients yet.
Still, as in the crystalline setting, this complex is endowed with a global Frobenius lift.

In rigid cohomology, two categories play an important role: convergent and overconvergent $F$-isocrystals.
Despite not having a $6$-functor formalism, they are a powerful tool; for instance, they are used to describe arithmetic $\mathcal D$-modules \cite{devissagedes}.

Inspired by the crystalline setting, a preprint of Ertl suggests that locally free overconvergent $F$-isocrystals can also be interpreted as coefficients for overconvergent de Rham--Witt cohomology \cite{comparisonbetweenrigidandoverconvergentcohomologywithcoefficients}.

In this paper, the third and last of a series, we back up her claim: there is an equivalence of categories between overconvergent $F$-isocrystals, and a category of locally projective overconvergent de Rham--Witt connections.
Our methods also give an equivalence of categories between convergent $F$-isocrystals and a category of locally projective de Rham--Witt connections.
This allows one to consider both these categories using only Zariski topology, and a global Frobenius
lift.\section*{Acknowledgments}This
work generalises the main result of my PhD thesis, so it is my pleasure to thank my advisor Daniel Caro.
I am also very grateful to all the members of my jury Andreas Langer, Tobias Schmidt, Christine Huyghe, Andrea Pulita and Jérôme Poineau.

The author is under ``Contrat EDGAR-CNRS no 277952 UMR 6623, financé par la région Bourgogne-Franche-Comté.''.
This work has been financed by ANR-21-CE39-0009-BARRACUDA.
It was also partially written during my position at Université Paris-Saclay, UVSQ, CNRS, Laboratoire de mathématiques de Versailles, 78000, Versailles, France.

Further acknowledgments to be added after the referee
process.\section{Connections}Throughout
this article, we let $p$ be a prime number, and $k$ be a perfect field of characteristic $p$.

In this section, $X$ shall denote a topological space, $R$ be a commutative ring, and $\mathcal F=\mleft(\bigoplus_{n\in\mathbb N}\mathcal F_n,d\mright)$ be a sheaf of strictly commutative $R$-dgas on $X$.
We first extend \cite[définition 2.2.1]{differentiellesnoncommutativesettheoriedegaloisdifferentielleouauxdifférences} to the context of
sheaves.\begin{enonce}{Definition}A
\emph{$\mathcal F$-connection} is a $\mathcal F_0$-module $M$ endowed with a morphism $\nabla\colon M\to M\otimes_{\mathcal F_0}\mathcal F_1$ of sheaves of groups over $X$.
We furthermore require that there exists an open cover $\mathcal C$ of $X$ such that for every open subset $U$ of an open in $\mathcal C$, the morphism $\nabla\mleft(U\mright)$ factors through a group morphism $\widetilde\nabla\colon M\mleft(U\mright)\to M\mleft(U\mright)\otimes_{\mathcal F_0\mleft(U\mright)}\mathcal F_1\mleft(U\mright)$ satisfying the Leibniz
rule:\begin{equation*}\forall s\in\mathcal F_0\mleft(U\mright),\ \forall m\in M\mleft(U\mright),\ \widetilde\nabla\mleft(sm\mright)=s\widetilde\nabla\mleft(m\mright)+m\otimes d\mleft(s\mright)\text.\end{equation*}\end{enonce}The
morphisms are defined as
usual.\begin{enonce}{Definition}Let
$\mleft(M,\nabla_M\mright)$ and $\mleft(N,\nabla_N\mright)$ two $\mathcal F$-connections.
A \emph{horizontal morphism} is a morphism $\varphi\colon M\to N$ of $\mathcal F_0$-modules such that the following diagram
commutes:\begin{equation*}\begin{tikzcd}[column sep=huge]M\arrow[r,"\varphi"]\arrow[d,"\nabla_M"]&N\arrow[d,"\nabla_N"]\\M\otimes_{\mathcal F_0}\mathcal F_1\arrow[r,"\varphi\otimes_{\mathcal F_0}\operatorname{Id}_{\mathcal F_1}"]&N\otimes_{\mathcal F_0}\mathcal F_1\text.\end{tikzcd}\end{equation*}

We shall denote by $\mathcal F\textsf-\mathsf{MC}$ the category of $\mathcal F$-connections and horizontal
morphisms.\end{enonce}We
we endow the category of graded abelian groups with the usual translation functor $\bigoplus_{i\in\mathbb N}G_i\mapsto\bigoplus_{i\in\mathbb N}G_{i+1}$.
We will denote by $\mathsf{dGAb}$ the category of differential objects in graded abelian groups.

We will construct a functor $\mathcal F\textsf-\mathsf{MC}\to\mathsf{Sh}\mleft(X,\mathsf{dGAb}\mright)$ from the category of $\mathcal F$-connections to the category of sheaves on $X$ of differential objects in graded abelian groups.
To define it, on any cover $\mathcal C$ satisfying the assumptions of the definition of a $\mathcal F$-connection, and on any open $U$ contained in an open of $\mathcal C$, for each $k\in\mathbb N$
put:\begin{equation*}\widetilde\nabla_k\colon\begin{array}{rl}M\mleft(U\mright)\otimes_{\mathcal F_0\mleft(U\mright)}\mathcal F_k\mleft(U\mright)\to&M\mleft(U\mright)\otimes_{\mathcal F_0\mleft(U\mright)}\mathcal F_{k+1}\mleft(U\mright)\\m\otimes s\mapsto&\widetilde\nabla\mleft(m\mright)\wedge s+m\otimes d\mleft(s\mright)\end{array}\end{equation*}

These maps are well-defined, by applying the same computations as in \cite[07I0]{stacksproject}.
Hence, for each open $V\in\mathcal C$ they yield glueable morphisms of sheaves of groups $\nabla_k|_V\colon\mleft(M\otimes_{\mathcal F_0}\mathcal F_k\mright)|_V\to\mleft(M\otimes_{\mathcal F_0}\mathcal F_{k+1}\mright)|_V$, which give rise to the sheaf $\mleft(M\otimes_{\mathcal F_0}\mathcal F,\nabla_\bullet\mright)$ of $\mathsf{Sh}\mleft(X,\mathsf{dGAb}\mright)$.
Of course,
$\nabla_0=\nabla$.\begin{enonce}{Definition}A
$\mathcal F$-connection $\nabla$ is said to be \emph{integrable} if $\nabla_1\circ\nabla_0=0$.
We shall denote by $\mathcal F\textsf-\mathsf{MIC}$ the full subcategory of $\mathcal F\textsf-\mathsf{MC}$ whose objects are integrable
connections.\end{enonce}The
computations in the proofs of \cite[proposition II 3.2.3 and proposition II 3.2.5]{cohomologiecristallinedesschemasdecaracteristique} still hold in this new context.
In particular, this means that a $\mathcal F$-connection $\nabla$ is integrable if and only if $\mleft(M\otimes_{\mathcal F_0}\mathcal F,\nabla_\bullet\mright)$ is a sheaf of complexes in graded abelian groups, thus if and only if it is a right differential graded
$\mathcal F$-module.\begin{enonce}{Proposition}The
categories $\mathcal F\textsf-\mathsf{MC}$ and $\mathcal F\textsf-\mathsf{MIC}$ are
additive.\end{enonce}\begin{proof}It
is straightforward to check that the addition of two horizontal morphisms is still horizontal, so the categories are preadditive.
Moreover, the product of a finite family of $\mathcal F_0$-modules endowed with (integrable) $\mathcal F$-connections is also endowed with a canonical (integrable) $\mathcal F$-connection since the tensor products and direct sums commute.
It is easily checked that the projections are horizontal, and that this new connection is the finite product in these
categories.\end{proof}\begin{enonce}{Definition}Let
$\nabla$ be a $\mathcal F$-connection on a free $\mathcal F_0$-module $M$.
Let $\mleft\{m_i\mright\}_{i\in I}$ be a $\mathcal F_0\mleft(X\mright)$-basis of $M\mleft(X\mright)$, where $I$ is a set.

The unique matrix $\mleft(n_{i,j}\mright)_{i,j\in I}\in\operatorname{Mat}_{\#I}\mleft(\mathcal F_1\mleft(X\mright)\mright)$ such that for every $i\in I$ we have $\nabla\mleft(m_i\mright)=\sum_{j\in I}m_j\otimes n_{j,i}$ is called the \emph{representative matrix} of $\nabla$ for the basis $\mleft\{m_i\mright\}_{i\in I}$.

Similarly, the unique matrix $\mleft(n_{i,j}\mright)_{i,j\in I}\in\operatorname{Mat}_{\#I}\mleft(\mathcal F_2\mleft(X\mright)\mright)$ such that for every $i\in I$ we have $\nabla_1\circ\nabla_0\mleft(m_i\mright)=\sum_{j\in I}m_j\otimes n_{j,i}$ is called the \emph{curvature matrix} of $\nabla$ for the basis
$\mleft\{m_i\mright\}_{i\in I}$.\end{enonce}A
representative matrix uniquely determines a $\mathcal F$-connection on a free $\mathcal F_0$-module for a given basis, and conversely using the following
result.\begin{enonce}{Proposition}\label{connectionevaluation}Let
$\nabla\colon M\to M\otimes_{\mathcal F_0}\mathcal F_1$ be a $\mathcal F$-connection on a free $\mathcal F_0$-module.
Let $\mleft\{m_i\mright\}_{i\in I}$ be a $\mathcal F_0\mleft(X\mright)$-basis of $M\mleft(X\mright)$, where $I$ is a set.
Let $U=\mleft\{u_i\mright\}_{i\in I}\in\mathcal F_0\mleft(X\mright)^I$ be a vertical vector, and let $N=\mleft(n_{i,j}\mright)_{i,j\in I}\in\operatorname{Mat}_{\#I}\mleft(\mathcal F_1\mleft(X\mright)\mright)$ be the representative matrix of $\nabla$ for $\mleft\{m_i\mright\}_{i\in I}$.

Then $\nabla\mleft(\sum_{i\in I}u_im_i\mright)=\sum_{i\in I}m_i\otimes v_i$,
where:\begin{equation*}\mleft\{v_i\mright\}_{i\in I}=NU+d\mleft(U\mright)\text.\end{equation*}\end{enonce}\begin{proof}Clear
from the Leibniz
rule.\end{proof}\begin{enonce}{Proposition}\label{connectionbasechange}Let
$\nabla\colon M\to M\otimes_{\mathcal F_0}\mathcal F_1$ be a $\mathcal F$-connection on a free $\mathcal F_0$-module.
Let $\mleft\{m_i\mright\}_{i\in I}$ be a $\mathcal F_0\mleft(X\mright)$-basis of $M\mleft(X\mright)$, where $I$ is a set.
Consider an invertible matrix $U=\mleft\{u_{i,j}\mright\}_{i,j\in I}\in\operatorname{GL}_{\#I}\mleft(\mathcal F_0\mleft(X\mright)\mright)$, and let $N=\mleft(n_{i,j}\mright)_{i,j\in I}\in\operatorname{Mat}_{\#I}\mleft(\mathcal F_1\mleft(X\mright)\mright)$ be the representative matrix of $\nabla$ for $\mleft\{m_i\mright\}_{i\in I}$.

The representative matrix of $\nabla$ for the basis $\mleft\{\sum_{j\in I}u_{j,i}m_j\mright\}_{i\in I}$ is
then:\begin{equation*}U^{-1}NU+U^{-1}d\mleft(U\mright)\text.\end{equation*}\end{enonce}\begin{proof}This
is a simple consequence of proposition
\ref{connectionevaluation}.\end{proof}\begin{enonce}{Proposition}\label{curvatureislinear}Let
$\nabla$ be a $\mathcal F$-connection.
Then, for each $k\in\mathbb N$ the morphism $\nabla_{k+1}\circ\nabla_k$ is
$\mathcal F_0$-linear.\end{enonce}\begin{proof}Just
follow the end of the proof of \cite[proposition II 3.2.3]{cohomologiecristallinedesschemasdecaracteristique}, and acknowledge that it still holds in our slightly different
setting.\end{proof}\begin{enonce}{Proposition}\label{integrabilitycharacterisation}Let
$\nabla\colon M\to M\otimes_{\mathcal F_0}\mathcal F_1$ be a $\mathcal F$-connection on a free $\mathcal F_0$-module.
Let $\mleft\{m_i\mright\}_{i\in I}$ be a $\mathcal F_0\mleft(X\mright)$-basis of $M\mleft(X\mright)$, where $I$ is a set, and let $N$ be a representative matrix of $\nabla$ for this basis.
Then, the curvature matrix of $\nabla$ for $\mleft\{m_i\mright\}_{i\in I}$
is:\begin{equation*}N^2+d\mleft(N\mright)\text.\end{equation*}

In particular, $\nabla$ is integrable if and only if $d\mleft(N\mright)=-N^2$.\end{enonce}\begin{proof}Write
$\mleft\{n_{i,j}\mright\}_{i,j\in I}\coloneqq N$.
By definition, $\nabla\mleft(m_i\mright)=\sum_{j\in I}m_j\otimes n_{j,i}$ for every $i\in I$.
Hence, we
find:\begin{equation*}\begin{split}\nabla_1\circ\nabla_0\mleft(m_i\mright)&=\nabla_1\mleft(\sum_{j\in I}m_j\otimes n_{j,i}\mright)\\&=\sum_{j\in I}\mleft(\nabla\mleft(m_j\mright)n_{j,i}+m_j\otimes d\mleft(n_{j,i}\mright)\mright)\\&=\sum_{j\in I}\mleft(\sum_{k\in I}m_k\otimes n_{k,j}n_{j,i}+m_j\otimes d\mleft(n_{j,i}\mright)\mright)\\\nabla_1\circ\nabla_0\mleft(m_i\mright)&=\sum_{j\in I}m_j\otimes\mleft(\sum_{k\in I}n_{j,k}n_{k,i}+d\mleft(n_{j,i}\mright)\mright)\text.\end{split}\end{equation*}

The second part of the statement is an application of proposition
\ref{curvatureislinear}.\end{proof}\begin{enonce}{Proposition}\label{functorfrommorphismconnections}Let
$\varphi\colon\mathcal F\to\mathcal G$ be a morphism of sheaves of strictly commutative $R$-dgas on $X$.
Let $\nabla\colon M\to M\otimes_{\mathcal F_0}\mathcal F_1$ be a $\mathcal F$-connection.
For every open $U$ of $X$,
let:\begin{equation*}\varphi_M\colon\begin{array}{rl}M\mleft(U\mright)\otimes_{\mathcal F_0\mleft(U\mright)}\mathcal F_1\mleft(U\mright)\to&\mleft(M\mleft(U\mright)\otimes_{\mathcal F_0\mleft(U\mright)}\mathcal G_0\mleft(U\mright)\mright)\otimes_{\mathcal G_0\mleft(U\mright)}\mathcal G_1\mleft(U\mright)\\m\otimes f\mapsto&\mleft(m\otimes1\mright)\otimes\varphi\mleft(f\mright)\end{array}\text.\end{equation*}

Then, assuming that $\nabla$ factors through $\widetilde\nabla\colon M\mleft(U\mright)\to M\mleft(U\mright)\otimes_{\mathcal F_0\mleft(U\mright)}\mathcal F_1\mleft(U\mright)$ satisfying the Leibniz rule on $U$, then the data of the following maps for varying $U$ covering $X$ yields a $\mathcal G$-connection $\varphi^*\mleft(\nabla\mright)$ on
$M\otimes_{\mathcal F_0}\mathcal G_0$:\begin{equation*}\widetilde{\varphi^*\mleft(\nabla\mright)}\mleft(U\mright)\colon\begin{array}{rl}M\mleft(U\mright)\otimes_{\mathcal F_0\mleft(U\mright)}\mathcal G_0\mleft(U\mright)\to&\mleft(M\mleft(U\mright)\otimes_{\mathcal F_0\mleft(U\mright)}\mathcal G_0\mleft(U\mright)\mright)\otimes_{\mathcal G_0\mleft(U\mright)}\mathcal G_1\mleft(U\mright)\\m\otimes g\mapsto&g\varphi_M\mleft(\widetilde\nabla\mleft(m\mright)\mright)+\mleft(m\otimes1\mright)\otimes d\mleft(g\mright)\end{array}\text.\end{equation*}

More precisely, it defines a
functor:\begin{equation*}\varphi^*\colon\mathcal F\textsf-\mathsf{MC}\to\mathcal G\textsf-\mathsf{MC}\text.\end{equation*}\end{enonce}\begin{proof}The
map is well-defined.
Indeed, it is biadditive and for every open $U$ of $X$, every $f\in\mathcal F_0\mleft(U\mright)$, every $m\in M\mleft(U\mright)$ and every $g\in\mathcal G_0\mleft(U\mright)$ one
has:\begin{equation*}\begin{split}\widetilde{\varphi^*\mleft(\nabla\mright)}\mleft(U\mright)\mleft(fm\otimes g\mright)&=g\varphi_M\mleft(\widetilde\nabla\mleft(fm\mright)\mright)+\mleft(fm\otimes1\mright)\otimes d\mleft(g\mright)\\&=g\varphi_M\mleft(f\widetilde\nabla\mleft(m\mright)+m\otimes d\mleft(f\mright)\mright)+\mleft(fm\otimes1\mright)\otimes d\mleft(g\mright)\\&=\varphi\mleft(f\mright)g\varphi_M\mleft(\widetilde\nabla\mleft(m\mright)\mright)+\mleft(m\otimes1\mright)\otimes\mleft(g\varphi\mleft(d\mleft(f\mright)\mright)+\varphi\mleft(f\mright)d\mleft(g\mright)\mright)\\\widetilde{\varphi^*\mleft(\nabla\mright)}\mleft(U\mright)\mleft(fm\otimes g\mright)&=\widetilde{\varphi^*\mleft(\nabla\mright)}\mleft(U\mright)\mleft(m\otimes\varphi\mleft(f\mright)g\mright)\text.\end{split}\end{equation*}

Moreover, we get a $\mathcal G$-connection because for all $m\in M$ and all $g,g'\in\mathcal G_0\mleft(U\mright)$ one
finds:\begin{equation*}\begin{split}\widetilde{\varphi^*\mleft(\nabla\mright)}\mleft(U\mright)\mleft(g'\mleft(m\otimes g\mright)\mright)&=\widetilde{\varphi^*\mleft(\nabla\mright)}\mleft(U\mright)\mleft(m\otimes g'g\mright)\\&=g'g\varphi_M\mleft(\widetilde\nabla\mleft(m\mright)\mright)+\mleft(m\otimes1\mright)\otimes d\mleft(g'g\mright)\\&=g'g\varphi_M\mleft(\widetilde\nabla\mleft(m\mright)\mright)+\mleft(m\otimes1\mright)\otimes g'd\mleft(g\mright)+\mleft(m\otimes1\mright)\otimes gd\mleft(g'\mright)\\\widetilde{\varphi^*\mleft(\nabla\mright)}\mleft(U\mright)\mleft(g'\mleft(m\otimes g\mright)\mright)&=g'\widetilde{\varphi^*\mleft(\nabla\mright)}\mleft(U\mright)\mleft(m\otimes g\mright)+\mleft(m\otimes g\mright)\otimes d\mleft(g'\mright)\text.\end{split}\end{equation*}

The functoriality is
clear.\end{proof}\begin{enonce}{Proposition}\label{functorfrommorphismintegrableconnections}Let
us keep the notations of proposition \ref{functorfrommorphismconnections}.
Then, the functor $\varphi^*$ preserves integrability; that is, by restriction we have a
functor:\begin{equation*}\varphi^*\colon\mathcal F\textsf-\mathsf{MIC}\to\mathcal G\textsf-\mathsf{MIC}\text.\end{equation*}\end{enonce}\begin{proof}By
construction, the functor $\mathcal F\textsf-\mathsf{MC}\to\mathsf{Sh}\mleft(X,\mathsf{dGAb}\mright)$ has the following property: when the $\mathcal F$-connection $\nabla\colon M\to M\otimes_{\mathcal F_0}\mathcal F_1$ is integrable, its image is the unique couple $\mleft(M\otimes_{\mathcal F_0}\mathcal F,\nabla_\bullet\mright)$ turning $M\otimes_{\mathcal F_0}\mathcal F$ into a right differential graded $\mathcal F$-module such that $\nabla_0=\nabla$.
Moreover, recall that $\mleft(M\otimes_{\mathcal F_0}\mathcal F,\nabla_\bullet\mright)$ is a right differential graded $\mathcal F$-module if and only if $\nabla$ is integrable.

Also, we have a tensor product on the category $\mathsf{RDiffGr}\textsf-\mathcal F\textsf-\mathsf{Mod}$ of right differential graded $\mathcal F$-modules \cite[09LL and 0FR2]{stacksproject} yielding a base change functor.
As can be seen from its definition, it coincides up to an isomorphism of graded $\mathcal F$-modules (truncated in degree greater than $1$) to the base change $\varphi^*$ we have defined.
Thus, the following diagram is essentially commutative (that is, commutative up to isomorphism of
functors):\begin{equation*}\begin{tikzcd}\mathcal F\textsf-\mathsf{MIC}\arrow[rr,"\varphi^*"]\arrow[d,"\nabla\mapsto\nabla_\bullet"]&&\mathcal G\textsf-\mathsf{MC}\arrow[d,"\nabla\mapsto\nabla_\bullet"]\\\mathsf{RDiffGr}\textsf-\mathcal F\textsf-\mathsf{Mod}\arrow[r,"\bullet\otimes_\mathcal F\mathcal G"]&\mathsf{RDiffGr}\textsf-\mathcal G\textsf-\mathsf{Mod}\arrow[r,"\mathsf{Forgetful}"]&\mathsf{Sh}\mleft(X,\mathsf{dGAb}\mright)\text.\end{tikzcd}\end{equation*}\end{proof}\section{Weak formal schemes}We
assume the reader is familiar with Meredith's theory of weak formal schemes, which he called weak formal preschemes in \cite[4. definition 1]{weakformalschemes}.
In what follows, the weak completion, in the sense of Monsky--Washnitzer \cite[definition 1.1]{formalcohomologyi}, of a commutative $W\mleft(k\mright)$-algebra $R$ shall be denoted by $R^\dagger$, and its ring of Witt vectors by $W\mleft(R\mright)$.
Its associated affine weak formal scheme shall be denoted by $\operatorname{Spff}\mleft(R^\dagger\mright)$.

In all the article, for every morphism of commutative rings $R\to S$ we shall
denote:\begin{equation*}\mathsf\Omega_{S/R}\coloneqq\Omega_{S/R}/\bigcap_{i\in\mathbb N}p^i\Omega_{S/R}\text.\end{equation*}

It is the universal $p$-adically separated de Rham complex associated to $R\to S$.
These similar notations shall lead to no confusion, as we shall always work with the $p$-adically separated complex, and never with the usual
one.\begin{enonce}{Proposition}\label{weakkahler}Let
$R$ be a commutative $W\mleft(k\mright)$-algebra such that $\mathsf\Omega^1_{R/W\mleft(k\mright)}$ is a finite $R$-module.
Then, we have a canonical isomorphism of
$R^\dagger$-modules:\begin{equation*}\mathsf\Omega^1_{R/W\mleft(k\mright)}\otimes_RR^\dagger\cong\mathsf\Omega^1_{R^\dagger/W\mleft(k\mright)}\text.\end{equation*}\end{enonce}\begin{proof}Let
$n\in\mathbb N$ and $\mleft(P_j\mright)_{j\in\mathbb N}\in\mleft(W\mleft(k\mright)\mleft[X_1,\ldots,X_n\mright]\mright)^\mathbb N$ such that there exists $c\in\mathbb N$ satisfying $\deg\mleft(P_j\mright)\leqslant c\mleft(j+1\mright)$ for every $j\in\mathbb N$.
If $\underline r\in R^n$, then by definition $\sum_{j\in\mathbb N}p^jP_j\mleft(\underline r\mright)$ is in $R^\dagger$, as every element of this ring can be written as such.
If $d\colon R\to\mathsf\Omega^1_{R/W\mleft(k\mright)}$ is the canonical derivation, then for every $j\in\mathbb N$ we
have:\begin{equation*}d\mleft(P_j\mleft(\underline r\mright)\mright)=\sum_{i=1}^nP_{j,i}\mleft(\underline r\mright)d\mleft(r_i\mright)\end{equation*}where
$P_{j,i}\in W\mleft(k\mright)\mleft[X_1,\ldots,X_n\mright]$ satisfies $\deg\mleft(P_{j,i}\mright)\leqslant c\mleft(j+1\mright)-1$ for every $i\in\mleft\llbracket1,n\mright\rrbracket$.
We
let:\begin{equation*}d'\colon\begin{array}{rl}R^\dagger\to&\mathsf\Omega^1_{R/W\mleft(k\mright)}\otimes_RR^\dagger\\\sum_{j\in\mathbb N}p^jP_j\mleft(\underline r\mright)\mapsto&\sum_{i=1}^nd\mleft(r_i\mright)\otimes\sum_{j\in\mathbb N}p^jP_{j,i}\mleft(\underline r\mright)\end{array}\end{equation*}

This map does not depend on neither $n$, nor $\mleft(P_j\mright)_{j\in\mathbb N}$ nor $\underline r$.

Indeed, $\mathsf\Omega^1_{R/W\mleft(k\mright)}\otimes_RR^\dagger$ is a finite $R^\dagger$-module.
So by \cite[theorem 1.6]{formalcohomologyi} and \cite[chap. III, §3, proposition 6]{algebrecommutativechapitres}, it is also $p$-adically separated.
Hence, if we were to choose other $n'\in\mathbb N$, $\mleft({P'}_j\mright)_{j\in\mathbb N}\in\mleft(W\mleft(k\mright)\mleft[X_1,\ldots,X_n\mright]\mright)^\mathbb N$ and $\underline r'\in R^n$ satisfying the same conditions, and such that $\sum_{j\in\mathbb N}p^jP_j\mleft(\underline r\mright)=\sum_{j\in\mathbb N}p^j{P'}_j\mleft(\underline r'\mright)$, then for every $k\in\mathbb N$ we
get:\begin{multline*}\sum_{i=1}^nd\mleft(r_i\mright)\otimes\sum_{j\in\mathbb N}p^jP_{j,i}\mleft(\underline r\mright)-\sum_{i=1}^{n'}d\mleft({r'}_i\mright)\otimes\sum_{j\in\mathbb N}p^j{P'}_{j,i}\mleft(\underline r'\mright)\\=\sum_{j=0}^{k-1}\mleft(\sum_{i=1}^np^jP_{j,i}\mleft(\underline r\mright)d\mleft(r_i\mright)-\sum_{i=1}^{n'}p^j{P'}_{j,i}\mleft(\underline r'\mright)d\mleft({r'}_i\mright)\mright)\otimes1\\+p^k\sum_{j=k}^{+\infty}\mleft(\sum_{i=1}^nd\mleft(r_i\mright)\otimes p^{j-k}P_{j,i}\mleft(\underline r\mright)-\sum_{i=1}^{n'}d\mleft({r'}_i\mright)\otimes p^{j-k}{P'}_{j,i}\mleft(\underline r'\mright)\mright)\text.\end{multline*}

Hence, the difference is divisible by $p^k$ for every $k\in\mathbb N$, so the map is well defined.

We see that $d'$ is a $W\mleft(k\mright)$-derivation because for every $j,j'\in\mathbb N$ we
have:\begin{equation*}\begin{split}d\mleft(P_j\mleft(\underline r\mright){P'}_{j'}\mleft(\underline r'\mright)\mright)&=P_j\mleft(\underline r\mright)d\mleft({P'}_{j'}\mleft(\underline r'\mright)\mright)+{P'}_{j'}\mleft(\underline r'\mright)d\mleft(P_j\mleft(\underline r\mright)\mright)\\&=\sum_{i=1}^{n'}P_j\mleft(\underline r\mright){P'}_{j',i}\mleft(\underline r'\mright)d\mleft({r'}_i\mright)+\sum_{i=1}^n{P'}_{j'}\mleft(\underline r'\mright)P_{j,i}\mleft(\underline r\mright)d\mleft(r_i\mright)\text.\end{split}\end{equation*}

If we let $d_\dagger\colon R^\dagger\to\mathsf\Omega^1_{R^\dagger/W\mleft(k\mright)}$ be the universal derivation, we then know there exists a unique morphism of $R^\dagger$-modules $\varphi\colon\mathsf\Omega^1_{R^\dagger/W\mleft(k\mright)}\to\mathsf\Omega^1_{R/W\mleft(k\mright)}\otimes_RR^\dagger$ such that $\varphi\circ d_\dagger=d'$.
Let $\iota\colon\mathsf\Omega^1_{R/W\mleft(k\mright)}\to\mathsf\Omega^1_{R^\dagger/W\mleft(k\mright)}$ be the canonical morphism of $R$-modules.
We
let:\begin{equation*}\varphi^{-1}\colon\begin{array}{rl}\mathsf\Omega^1_{R/W\mleft(k\mright)}\otimes_RR^\dagger\to&\mathsf\Omega^1_{R^\dagger/W\mleft(k\mright)}\\\omega\otimes r\mapsto&\iota(\omega)r\end{array}\end{equation*}

For every $r\in R$ we have $d'\mleft(r\mright)=d\mleft(r\mright)\otimes1$, which yields $\varphi\circ\varphi^{-1}\circ d'\mleft(r\mright)=d'\mleft(r\mright)$ and $\varphi^{-1}\circ\varphi\circ d_\dagger\mleft(r\mright)=d_\dagger\mleft(r\mright)$.
By using $p$-adically overconvergent series, one gets similar equalities when $r\in R^\dagger$.
Since $d'\mleft(R^\dagger\mright)$ is a $R^\dagger$-generating set of $\mathsf\Omega^1_{R/W\mleft(k\mright)}\otimes_RR^\dagger$ and $d_\dagger\mleft(R^\dagger\mright)$ is a $R^\dagger$-generating set of $\mathsf\Omega^1_{R^\dagger/W\mleft(k\mright)}$, one can
conclude.\end{proof}\begin{enonce}{Proposition}\label{derhamweak}Let
$k$ be a perfect field of characteristic $p$.
Assume that $R$ is a commutative $W\mleft(k\mright)$-algebra such that $\mathsf\Omega^1_{R/W\mleft(k\mright)}$ is a finite $R$-module.
Then, there is a canonical isomorphism of
$W\mleft(k\mright)$-dgas:\begin{equation*}\mathsf\Omega_{R/W\mleft(k\mright)}\otimes_RR^\dagger\cong\mathsf\Omega_{R^\dagger/W\mleft(k\mright)}\text.\end{equation*}\end{enonce}\begin{proof}This
is proposition \ref{weakkahler} combined with
\cite[III.83, proposition 8]{algebrechapitres}.\end{proof}The
convergent version of this statement, that is $\mathsf\Omega_{R/W\mleft(k\mright)}\otimes_R\widehat R\cong\mathsf\Omega_{\widehat R/W\mleft(k\mright)}$ with the above assumptions, is well-known by \cite[00RV and 0315]{stacksproject} and \cite[chap. III, §3, proposition 6]{algebrecommutativechapitres}.

We recall the following definition of Caro
\cite[définition 1.2.1]{devissagedes}.\begin{enonce}{Definition}A
weak formal scheme $\mathfrak X^\dagger$ over $W\mleft(k\mright)$ is said to be \emph{smooth} if for every $i\in\mathbb N$, the scheme $\mleft(\mathfrak X^\dagger,\mathcal O_{\mathfrak X^\dagger}/p^i\mathcal O_{\mathfrak X^\dagger}\mright)$ is smooth over $W\mleft(k\mright)/p^iW\mleft(k\mright)$.\end{enonce}

If $\mleft(X,\mathcal O_X\mright)$ is a ringed space, for every $\mathcal O_X\mleft(X\mright)$-module $M$ we denote by $\widetilde M$ the sheaf of $\mathcal O_X$-modules associated to $M$.
We refer to \cite[01BH]{stacksproject} for its properties.
This sheaf depends on the choice of the ringed space, but the notation should lead to no confusion in context.

If $\mleft(\mathfrak X^\dagger,\mathcal O_{\mathfrak X^\dagger}\mright)$ is an affine weak formal scheme, then by \cite[3. theorem 3]{weakformalschemes} we know that the functor $\widetilde\bullet$ induces an equivalence of categories between $\mathcal O_{\mathfrak X^\dagger}\mleft(\mathfrak X^\dagger\mright)$-modules of finite type and coherent $\mathcal O_{\mathfrak X^\dagger}$-modules.

The following proposition explains us that, in the context of weak affine schemes, one can extend this equivalence of categories to see that our notion of connection coincides with the usual algebraic one; so we will identify both point of views in the
paper.\begin{enonce}{Proposition}Let
$k$ be a perfect field of characteristic $p$.
Let $\mleft(\mathfrak X^\dagger,\mathcal O_{\mathfrak X^\dagger}\mright)$ be a smooth affine weak formal scheme over $W\mleft(k\mright)$.
Then, there exists a unique quasi-coherent sheaf of strictly commutative $W\mleft(k\mright)$-dgas $\mathsf\Omega_{\mathfrak X^\dagger/W\mleft(k\mright)}$ on $\mathfrak X^\dagger$ such that for every $\mathcal O_{\mathfrak X^\dagger}\mleft(\mathfrak X^\dagger\mright)$-module $M$ and every affine open set $U\subset\mathfrak X^\dagger$ we
have:\begin{equation*}\widetilde M\otimes_{\mathcal O_{\mathfrak X^\dagger}}\mathsf\Omega_{\mathfrak X^\dagger/W\mleft(k\mright)}\mleft(U\mright)\cong M\otimes_{\mathcal O_{\mathfrak X^\dagger}\mleft(\mathfrak X^\dagger\mright)}\mathsf\Omega_{\mathcal O_{\mathfrak X^\dagger}\mleft(U\mright)/W\mleft(k\mright)}\text.\end{equation*}\end{enonce}\begin{proof}Let
$U$ be an affine open set of $\mathfrak X^\dagger$.
Let $M$ be a $\mathcal O_{\mathfrak X^\dagger}\mleft(\mathfrak X^\dagger\mright)$-module.
Let $\mathsf\Omega_{\mathfrak X^\dagger/W\mleft(k\mright)}$ be the sheafification of $\mathsf\Omega_{\mathcal O_{\mathfrak X^\dagger}\mleft(\bullet\mright)/W\mleft(k\mright)}$.

The canonical morphism $M\otimes_{\mathcal O_{\mathfrak X^\dagger}\mleft(\mathfrak X^\dagger\mright)}\mathsf\Omega_{\mathcal O_{\mathfrak X^\dagger}\mleft(U\mright)/W\mleft(k\mright)}\to\widetilde M\otimes_{\mathcal O_{\mathfrak X^\dagger}}\mathcal\mathsf\Omega_{\mathfrak X^\dagger/W\mleft(k\mright)}\mleft(U\mright)$ of $\mathcal O_{\mathfrak X^\dagger}\mleft(U\mright)$-modules yields by adjunction \cite[01BH]{stacksproject} a morphism of sheaves of $\mathcal O_{\mathfrak X^\dagger}|_U$-modules $\widetilde{M\otimes_{\mathcal O_{\mathfrak X^\dagger}\mleft(\mathfrak X^\dagger\mright)}\mathsf\Omega_{\mathcal O_{\mathfrak X^\dagger}\mleft(U\mright)/W\mleft(k\mright)}}\to\widetilde M\otimes_{\mathcal O_{\mathfrak X^\dagger}}\mathsf\Omega_{\mathfrak X^\dagger/W\mleft(k\mright)}|_U$.
For every $x\in U$, this morphism induces a morphism of
$\mathcal O_{\mathfrak X^\dagger,x}$-modules:\begin{multline}\label{kahlerongerms}M\otimes_{\mathcal O_{\mathfrak X^\dagger}\mleft(\mathfrak X^\dagger\mright)}\mathsf\Omega_{\mathcal O_{\mathfrak X^\dagger}\mleft(U\mright)/W\mleft(k\mright)}\otimes_{\mathcal O_{\mathfrak X^\dagger}\mleft(U\mright)}\operatorname{colim}_{x\in V\subset U}\mathcal O_{\mathfrak X^\dagger}\mleft(V\mright)\\\to M\otimes_{\mathcal O_{\mathfrak X^\dagger}\mleft(\mathfrak X^\dagger\mright)}\operatorname{colim}_{x\in V\subset U}\mathsf\Omega_{\mathcal O_{\mathfrak X^\dagger}\mleft(V\mright)/W\mleft(k\mright)}\text.\end{multline}

One can restrict the above open sets $V$ to standard opens $V=D\mleft(\overline f\mright)$, where $f\in\mathcal O_{\mathfrak X^\dagger}\mleft(U\mright)$ satisfies $x\in D\mleft(\overline f\mright)$.
By \cite[00RT]{stacksproject} and \cite[III.83, proposition 8]{algebrechapitres}, there is an isomorphism $\Omega_{\mathcal O_{\mathfrak X^\dagger}\mleft(U\mright)/W\mleft(k\mright)}\otimes_{\mathcal O_{\mathfrak X^\dagger}\mleft(U\mright)}\mathcal O_{\mathfrak X^\dagger}\mleft(U\mright)_f\cong\Omega_{\mathcal O_{\mathfrak X^\dagger}\mleft(U\mright)_f/W\mleft(k\mright)}$ of $W\mleft(k\mright)$-dgas.
By applying proposition \ref{derhamweak}, one finds an isomorphism of
$W\mleft(k\mright)$-dgas:\begin{equation*}\mathsf\Omega_{\mathcal O_{\mathfrak X^\dagger}\mleft(U\mright)/W\mleft(k\mright)}\otimes_{\mathcal O_{\mathfrak X^\dagger}\mleft(U\mright)}\mathcal O_{\mathfrak X^\dagger}\mleft(V\mright)\cong\mathsf\Omega_{\mathcal O_{\mathfrak X^\dagger}\mleft(V\mright)/W\mleft(k\mright)}\text.\end{equation*}

Moreover, as tensor products commute with direct limits we see that for every $x\in U$, the morphism \eqref{kahlerongerms} is an isomorphism.
In particular, we have an isomorphism $\widetilde{M\otimes_{\mathcal O_{\mathfrak X^\dagger}\mleft(\mathfrak X^\dagger\mright)}\mathsf\Omega_{\mathcal O_{\mathfrak X^\dagger}\mleft(U\mright)/W\mleft(k\mright)}}\to\widetilde M\otimes_{\mathcal O_{\mathfrak X^\dagger}}\mathsf\Omega_{\mathfrak X^\dagger/W\mleft(k\mright)}|_U$ and we
conclude.\end{proof}Of
course, there is a similar statement for formal schemes with the usual $p$-adic
completion.\section{The overconvergent de Rham--Witt complex}We
let $\overline R$ be a commutative $k$-algebra of finite type.
In this paper, we shall denote by $W\Omega_{\overline R/k}$ the de Rham--Witt complex of $\overline R$.
We refer to \cite{complexedederhamwittetcohomologiecristalline} and \cite{derhamwittcohomologyforaproperandsmoothmorphism} for the definition and basic properties.

Let $n\in\mathbb N$.
In this article, we shall denote $k\mleft[\underline X\mright]\coloneqq k\mleft[X_1,\ldots,X_n\mright]$, and similarly for $W\mleft(k\mright)\mleft[\underline X\mright]$.
The truncated de Rham--Witt complex of $\overline R$ will be denoted by $W_n\Omega_{\overline R/k}$.

Recall that when $\overline R$ is smooth, a theorem of Elkik tells us that we can find a smooth commutative $W\mleft(k\mright)$-algebra $R$ lifting $\overline R$ \cite[07M8]{stacksproject}.
Moreover, its weak completion $R^\dagger$ does not depend, up to an isomorphism of $W\mleft(k\mright)$-algebras, on the choice of $R$, and there is a morphism of rings $F\colon R^\dagger\to R^\dagger$ lifting the Frobenius and compatible with $F\colon W\mleft(k\mright)\to W\mleft(k\mright)$ \cite[theorem 2.4.4]{thecohomologyofmonskyandwashnitzer}.
This yields in turn a morphism of $\delta$-rings $t_F\colon R^\dagger\to W\mleft(\overline R\mright)$ by \cite[théorème 4]{joyal}, which extends by the universal property of the $p$-adically separated de Rham complex \cite[theorem 4.2]{formalcohomologyi} to a morphism of
$W\mleft(k\mright)$-dgas:\begin{equation*}t_{F}\colon\mathsf\Omega_{R^\dagger/W\mleft(k\mright)}\to W\Omega_{\overline R/k}\text.\end{equation*}

For every $\varepsilon>0$, there is a
map:\begin{equation*}\zeta_\varepsilon\colon W\Omega_{k\mleft[\underline X\mright]/k}\to\mathbb R\cup\mleft\{+\infty,-\infty\mright\}\text.\end{equation*}We
do not recall their technical definition yet; we rather need their properties which we state below.

Let $\varphi\colon k\mleft[\underline X\mright]\to\overline R$ be a surjective map of $k$-algebras from a polynomial ring.
This yields a morphism of $W\mleft(k\mright)$-dgas $\varphi^*\colon W\Omega_{k\mleft[\underline X\mright]/k}\to W\Omega_{\overline R/k}$.
We
put:\begin{equation*}\zeta_{\varepsilon,\varphi}\colon\begin{array}{rl}W\Omega_{\overline R/k}\to&\mathbb R\cup\mleft\{+\infty,-\infty\mright\}\\w\mapsto&\sup\mleft\{\zeta_\varepsilon\mleft(x\mright)\mid\varphi^*\mleft(x\mright)=w\mright\}\end{array}\text.\end{equation*}

First of all, by \cite[proposition 5.14]{localstructureoftheoverconvergentderhamwittcomplex} and the subsequent discussion, we can make the following definition which coincides with \cite[definition 1.1]{overconvergentderhamwittcohomology} in the case where $k$ is a perfect field, and which does not depend on the choice of
$\varphi$.\begin{enonce}{Definition}The
\emph{overconvergent de Rham--Witt complex of $\overline R$} is:\begin{equation*}W^\dagger\Omega_{\overline R/k}\coloneqq\mleft\{w\in W\Omega_{\overline R/k}\mid\exists\varepsilon>0,\ \zeta_{\varepsilon,\varphi}\mleft(w\mright)\neq-\infty\mright\}\text.\end{equation*}\end{enonce}

Recall that when $\overline R=k\mleft[\underline X\mright]$ we have a decomposition of the overconvergent de Rham--Witt complex as three graded sub-$W\mleft(k\mright)$-modules
\cite[(8)]{pseudovaluationsonthederhamwittcomplex}:\begin{equation}\label{derhamwittdecomposition}W^\dagger\Omega_{k\mleft[\underline X\mright]/k}=W^\dagger\Omega^\mathrm{int}_{k\mleft[\underline X\mright]/k}\oplus W^\dagger\Omega^\mathrm{frp}_{k\mleft[\underline X\mright]/k}\oplus d\mleft(W^\dagger\Omega^\mathrm{frp}_{k\mleft[\underline X\mright]/k}\mright)\text.\end{equation}

By \cite[3.3.2]{relevementdesalgebreslissesetdeleursmorphismes}, there exists a lift of the Frobenius on ${W\mleft(k\mright)\mleft[\underline X\mright]}^\dagger$, which we shall also denote by $F$, such that we can find a $\delta$-ring morphism ${W\mleft(k\mright)\mleft[\underline X\mright]}^\dagger\to R^\dagger$ lifting $\varphi$.
More details can be found in \cite[lemma 7.1]{localstructureoftheoverconvergentderhamwittcomplex}.
Moreover, there is a canonical isomorphism of $W\mleft(k\mright)$-dgas $\mathsf\Omega_{{W\mleft(k\mright)\mleft[\underline X\mright]}^\dagger/W\mleft(k\mright)}\cong W^\dagger\Omega^\mathrm{int}_{k\mleft[\underline X\mright]/k}$ \cite[proposition 7.11]{localstructureoftheoverconvergentderhamwittcomplex}, so that we shall identify them in what follows.

For $t\in\mathbb N$, there exist subsets $H\mleft(t\mright)\subset W^\dagger\Omega^{\mathrm{int},t}_{k\mleft[\underline X\mright]/k}$ and $G\mleft(t\mright)\subset W^\dagger\Omega^{\mathrm{frp},t}_{k\mleft[\underline X\mright]/k}$, whose definition we will recall when needed, allowing us to introduce the following
${W\mleft(k\mright)\mleft[\underline X\mright]}^\dagger$-modules:\begin{align*}W^\dagger\Omega^{\underline X-\mathrm{int},t}_{k\mleft[\underline X\mright]/k}&\coloneqq\mleft\{\sum_{e\in H\mleft(t\mright)}t_F\mleft(s_ee\mright)\mid\begin{array}l\exists C\in\mathbb R,\ \exists\varepsilon>0,\ \forall e\in H\mleft(t\mright),\\\mleft(s_e\in W^\dagger\Omega^{\mathrm{int},0}_{k\mleft[\underline X\mright]/k}\mright)\wedge\mleft(\zeta_\varepsilon\mleft(s_e\mright)+\zeta_\varepsilon\mleft(e\mright)\geqslant C\mright)\end{array}\mright\}\text,\\W^\dagger\Omega^{\underline X-\mathrm{frp},t}_{k\mleft[\underline X\mright]/k}&\coloneqq\mleft\{\sum_{e\in G\mleft(t\mright)}t_F\mleft(s_e\mright)e\mid\begin{array}l\exists C\in\mathbb R,\ \exists\varepsilon>0,\ \forall e\in G\mleft(t\mright),\\\mleft(s_e\in W^\dagger\Omega^{\mathrm{int},0}_{k\mleft[\underline X\mright]/k}\mright)\wedge\mleft(\zeta_\varepsilon\mleft(s_e\mright)+\zeta_\varepsilon\mleft(e\mright)\geqslant C\mright)\end{array}\mright\}\text.\end{align*}

And we
put:\begin{equation*}
W^\dagger\Omega^{\underline X}_{k\mleft[\underline X\mright]/k}\coloneqq W^\dagger\Omega^{\underline X-\mathrm{int}}_{k\mleft[\underline X\mright]/k}+W^\dagger\Omega^{\underline X-\mathrm{frp}}_{k\mleft[\underline X\mright]/k}+d\mleft(W^\dagger\Omega^{\underline X-\mathrm{frp}}_{k\mleft[\underline X\mright]/k}\mright)\text.
\end{equation*}

For details on these constructions, see \cite[9]{localstructureoftheoverconvergentderhamwittcomplex}.
Note that in that paper, there are hypothesis on the morphism $\varphi$, but we can safely assume here that they are satisfied.

By \cite[theorem 9.7]{localstructureoftheoverconvergentderhamwittcomplex}, the previous decomposition as graded sub-$W\mleft(k\mright)$-modules actually extends for general
$\overline R$:\begin{equation*}W^\dagger\Omega_{\overline R/k}=W^\dagger\Omega^\mathrm{int}_{\overline R/k}\oplus W^\dagger\Omega^\mathrm{frp}_{\overline R/k}\oplus d\mleft(W^\dagger\Omega^\mathrm{frp}_{\overline R/k}\mright)\text.\end{equation*}

More precisely, this decomposition depends on $\varphi$ as $W^\dagger\Omega^\mathrm{int}_{\overline R/k}=\varphi^*\mleft(W^\dagger\Omega^\mathrm{int}_{k\mleft[\underline X\mright]/k}\mright)$ and $W^\dagger\Omega^\mathrm{frp}_{\overline R/k}=\varphi^*\mleft(W^\dagger\Omega^\mathrm{frp}_{k\mleft[\underline X\mright]/k}\mright)$.
In practice, there shall be no confusion so we omit it in our notations.
This enables us to define the following
map:\begin{equation*}\check\zeta_{\varepsilon,\varphi}\colon\begin{array}{rl}W\Omega_{\overline R/k}\to&\mathbb R\cup\mleft\{+\infty,-\infty\mright\}\\w\mapsto&\sup\mleft\{\zeta_\varepsilon\mleft(x\mright)\mid x\in{\varphi^*}^{-1}\mleft(w\mright)\cap W^\dagger\Omega^{\underline X}_{k\mleft[\underline X\mright]/k}\mright\}\end{array}\text.\end{equation*}

For $x\in W^\dagger\Omega_{\overline R/k}$, we will denote by $x|_\mathrm{int}$ its projection to $W^\dagger\Omega^\mathrm{int}_{\overline R/k}$, by $x|_\mathrm{frp}$ its projection to $W^\dagger\Omega^\mathrm{frp}_{\overline R/k}$, by $x|_\mathrm{d(frp)}$ its projection to $d\mleft(W^\dagger\Omega^\mathrm{frp}_{\overline R/k}\mright)$, and by $x|_\mathrm{frac}$ its projection to $W^\dagger\Omega^\mathrm{frac}_{\overline R/k}\coloneqq W^\dagger\Omega^\mathrm{frp}_{\overline R/k}+d\mleft(W^\dagger\Omega^\mathrm{frp}_{\overline R/k}\mright)$.

A similar decomposition exists for the usual de Rham--Witt complex \cite[theorem 4.5]{localstructureoftheoverconvergentderhamwittcomplex}.
The construction work the same way as here, except that we use the $p$-adic completion instead of the weak completion.
We omit the details.
We shall keep the same notations in that context.

Applying \cite[proposition 7.14 and theorem 9.7]{localstructureoftheoverconvergentderhamwittcomplex}, we get $\delta>0$ such that for every $\varepsilon\in\mleft]0,\delta\mright]$ we
have:\begin{align}\forall x\in W^\dagger\Omega^{\mathrm{int},0}_{k\mleft[\underline X\mright]/k},\ \zeta_\varepsilon\mleft(t_F\mleft(x\mright)-x\mright)&\geqslant\zeta_\varepsilon\mleft(t_F\mleft(x\mright)\mright)+\frac34\text,\label{tfandvf}\\\forall x\in W^\dagger\Omega_{\overline R/k},\ \check\zeta_{\varepsilon,\varphi}\mleft(x|_\mathrm{int}\mright)&\geqslant\zeta_{\varepsilon,\varphi}\mleft(x\mright)-\frac12\text,\label{costofdecompositionint}\\\forall x\in W^\dagger\Omega_{\overline R/k},\ \check\zeta_{\varepsilon,\varphi}\mleft(x|_\mathrm{frp}\mright)&\geqslant\zeta_{\varepsilon,\varphi}\mleft(x\mright)-\frac12\text,\label{costofdecompositionfrp}\\\forall x\in W^\dagger\Omega_{\overline R/k},\ \check\zeta_{\varepsilon,\varphi}\mleft(x|_\mathrm{d(frp)}\mright)&\geqslant\zeta_{\varepsilon,\varphi}\mleft(x\mright)-\frac12\text,\label{costofdecompositiondfrp}\\\forall x\in W^\dagger\Omega^\mathrm{frp}_{\overline R/k},\ \check\zeta_{\varepsilon,\varphi}\mleft(x\mright)&=\check\zeta_{\varepsilon,\varphi}\mleft(d\mleft(x\mright)\mright)\text,\label{frpfromdfrp}\\\forall x\in W^\dagger\Omega_{\overline R/k},\ \zeta_{\varepsilon,\varphi}\mleft(x\mright)&\geqslant\check\zeta_{\varepsilon,\varphi}\mleft(x\mright)\text.\label{zetabiggerthancheckzeta}\end{align}

Also, by \cite[(9)]{pseudovaluationsonthederhamwittcomplex}, we
have:\begin{equation}\label{zetaepsilonandd}\forall\varepsilon>0,\ \forall x\in W^\dagger\Omega_{k\mleft[\underline X\mright]/k},\ \zeta_\varepsilon\mleft(d\mleft(x\mright)\mright)\geqslant\zeta_\varepsilon\mleft(x\mright)\text.\end{equation}\begin{enonce}{Proposition}\label{xfrpminoration}Let
$t\in\mathbb N$.
For every $\varepsilon\in\mleft]0,\delta\mright]$ and every $x\in W^\dagger\Omega^{\underline X-\mathrm{frp},t}_{k\mleft[\underline X\mright]/k}$ we have for the decomposition
\eqref{derhamwittdecomposition}:\begin{align*}\zeta_\varepsilon\mleft(x|_\mathrm{int}\mright)&\geqslant\zeta_\varepsilon\mleft(x\mright)+\frac34\text,\\\zeta_\varepsilon\mleft(x|_\mathrm{frp}\mright)&\geqslant\zeta_\varepsilon\mleft(x\mright)\text,\\\zeta_\varepsilon\mleft(x|_\mathrm{d(frp)}\mright)&\geqslant\zeta_\varepsilon\mleft(x\mright)+\frac34\text.\end{align*}\end{enonce}\begin{proof}For
every $e\in G\mleft(t\mright)$, we let $s_e\in W^\dagger\Omega^{\mathrm{int},0}_{k\mleft[\underline X\mright]/k}$ such
that:\begin{align*}\sum_{e\in G\mleft(t\mright)}t_F\mleft(s_e\mright)e&=x\text,\\\zeta_\varepsilon\mleft(s_e\mright)+\zeta_\varepsilon\mleft(e\mright)&\geqslant\zeta_\varepsilon\mleft(x\mright)\text.\end{align*}

Now, write $x=\sum_{e\in G\mleft(t\mright)}\mleft(t_F\mleft(s_e\mright)-s_e\mright)e+\sum_{e\in G\mleft(t\mright)}s_ee$.
Following the remark after \cite[proposition 9.1]{localstructureoftheoverconvergentderhamwittcomplex}, the second series is in $W^\dagger\Omega^{\mathrm{frp},t}_{k\mleft[\underline X\mright]/k}$.
So now we can conclude thanks to \eqref{tfandvf} and to the fact that $\zeta_\varepsilon$ is a pseudovaluation \cite[theorem
3.17]{pseudovaluationsonthederhamwittcomplex}.\end{proof}Given
$t,s\in\mathbb N$, we write $\mleft\llbracket t,s\mright\rrbracket\coloneqq\mleft[t,s\mright]\cap\mathbb N$.
We recall that a rng is defined like a ring, without the need of an identity for the
multiplication.\begin{enonce}{Lemma}\label{rngrelation}Let
$R$ be a rng.
Let $t\in\mathbb N$.
Let $\mleft(x_i\mright)_{i\in\mleft\llbracket0,t\mright\rrbracket},\mleft(y_i\mright)_{i\in\mleft\llbracket0,t\mright\rrbracket}\in R^{t+1}$.
Let $\mathcal F_{\mleft(\mleft(x_i,y_i\mright)\mright)_{i\in\mleft\llbracket0,t\mright\rrbracket}}$ denote the set of functions $\mleft\llbracket0,t\mright\rrbracket\to R$ such that for every $i\in\mleft\llbracket0,t\mright\rrbracket$ we have $f_i\in\mleft\{x_i+y_i,-y_i\mright\}$, and $f_i=-y_i$ for at least one such $i$.
Then:\begin{equation*}\prod_{i=0}^t\mleft(x_i+y_i\mright)=\prod_{i=0}^tx_i-\sum_{f\in\mathcal F_{\mleft(\mleft(x_i,y_i\mright)\mright)_{i\in\mleft\llbracket0,t\mright\rrbracket}}}\prod_{i=0}^tf_i\text.\end{equation*}\end{enonce}\begin{proof}This
is done by induction on $t\in\mathbb N$, the case $t=0$ being easy.
Assuming the result for some $t\in\mathbb N$, we
get:\begin{equation*}\prod_{i=0}^tx_i\mleft(x_{t+1}+y_{t+1}\mright)=\prod_{i=0}^{t+1}x_i-\mleft(\prod_{i=0}^t\mleft(x_i+y_i\mright)+\sum_{f\in\mathcal F_{\mleft(\mleft(x_i,y_i\mright)\mright)_{i\in\mleft\llbracket0,t\mright\rrbracket}}}\prod_{i=0}^tf_i\mright)\mleft(-y_{t+1}\mright)\text.\end{equation*}\end{proof}\begin{enonce}{Proposition}\label{tfandvfhigherdegree}Let
$t\in\mathbb N$.
Let $x\in W^\dagger\Omega^{\mathrm{int},t}_{k\mleft[\underline X\mright]/k}$.
For the decomposition \eqref{derhamwittdecomposition}, for every $\varepsilon\in\mleft]0,\delta\mright]$ we
have:\begin{align*}\zeta_\varepsilon\mleft(t_F\mleft(x\mright)|_\mathrm{int}\mright)&\geqslant\zeta_\varepsilon\mleft(t_F\mleft(x\mright)\mright)\text,\\\zeta_\varepsilon\mleft(t_F\mleft(x\mright)|_\mathrm{frp}\mright)&\geqslant\zeta_\varepsilon\mleft(t_F\mleft(x\mright)\mright)+\frac34\text,\\\zeta_\varepsilon\mleft(t_F\mleft(x\mright)|_\mathrm{d(frp)}\mright)&\geqslant\zeta_\varepsilon\mleft(t_F\mleft(x\mright)\mright)+\frac34\text.\end{align*}\end{enonce}\begin{proof}We
derive from \cite[proposition 9.1]{localstructureoftheoverconvergentderhamwittcomplex} that for every $e\in H\mleft(t\mright)$ there exists a unique $s_e\in W^\dagger\Omega^{\mathrm{int},0}_{k\mleft[\underline X\mright]/k}$ such
that:\begin{align*}\zeta_\varepsilon\mleft(s_e\mright)+\zeta_\varepsilon\mleft(e\mright)&\geqslant\zeta_\varepsilon\mleft(t_F\mleft(x\mright)\mright)\text,\\\sum_{e\in H\mleft(t\mright)}t_F\mleft(s_ee\mright)&=t_F\mleft(x\mright)\text.\end{align*}

By definition \cite[9]{localstructureoftheoverconvergentderhamwittcomplex}, each $e\in H\mleft(t\mright)$ is of the form $\prod_{i=1}^td\mleft(y_{e,i}\mright)$ where the $y_{e,i}\in W^\dagger\Omega^{\mathrm{int},0}_{k\mleft[\underline X\mright]/k}$ for each $i\in\mleft\llbracket1,n\mright\rrbracket$, and $H\mleft(t\mright)$ is a finite set.

Now, fix $e\in H\mleft(t\mright)$.
Let $x_0\coloneqq s_e$, $x_i\coloneqq d\mleft(y_{e,i}\mright)$ for $i\in\mleft\llbracket1,t\mright\rrbracket$, $y_0\coloneqq t_F\mleft(s_e\mright)-s_e$ and $y_i\coloneqq d\mleft(t_F\mleft(y_{e,i}\mright)-y_{e,i}\mright)$ for $i\in\mleft\llbracket1,t\mright\rrbracket$.
Using lemma \ref{rngrelation} and its notations we
get:\begin{equation*}t_F\mleft(s_ee\mright)=s_ee-\sum_{f\in\mathcal F_{\mleft(\mleft(x_i,y_i\mright)\mright)_{i\in\mleft\llbracket0,t\mright\rrbracket}}}\prod_{i=0}^tf_i\text.\end{equation*}

For each $f$ in the above sum, we have $\zeta_\varepsilon\mleft(\prod_{i=0}^tf_i\mright)\geqslant\zeta_\varepsilon\mleft(s_e\mright)+\zeta_\varepsilon\mleft(e\mright)+\frac34$ by formulas \eqref{tfandvf} and \eqref{zetaepsilonandd}, and since $\zeta_\varepsilon$ is a pseudovaluation \cite[theorem 3.17]{pseudovaluationsonthederhamwittcomplex}.
This enables us to conclude with the same argument because $s_ee$ is in the $W\mleft(k\mright)$-dga
$W^\dagger\Omega^{\mathrm{int},t}_{k\mleft[\underline X\mright]/k}$, and because $\zeta_\varepsilon$ is compatible with this decomposition \cite[definition 3.4]{pseudovaluationsonthederhamwittcomplex}.\end{proof}\begin{enonce}{Proposition}\label{controllingintfrp}Let
$t\in\mathbb N^*$ and $s\in\mathbb N$.
Let $x\in W^\dagger\Omega^{\mathrm{int},t}_{\overline R/k}$ and $y\in W^\dagger\Omega^{\mathrm{frp},s}_{\overline R/k}$.
For every $\varepsilon\in\mleft]0,\delta\mright]$ we
have:\begin{equation*}\check\zeta_{\varepsilon,\varphi}\mleft(xy\mright)\geqslant\check\zeta_{\varepsilon,\varphi}\mleft(x\mright)+\check\zeta_{\varepsilon,\varphi}\mleft(y\mright)+\frac14\text.\end{equation*}\end{enonce}\begin{proof}Let
$\mu>0$.
Let $\widetilde x\in W^\dagger\Omega^{\underline X-\mathrm{int},t}_{k\mleft[\underline X\mright]/k}$ and $\widetilde y\in W^\dagger\Omega^{\underline X-\mathrm{frp},s}_{k\mleft[\underline X\mright]/k}$ be respective preimages of $x$ and $y$ by $\varphi^*$ such that $\zeta_\varepsilon\mleft(\widetilde x\mright)\geqslant\check\zeta_{\varepsilon,\varphi}\mleft(x\mright)-\mu$ and $\zeta_\varepsilon\mleft(\widetilde y\mright)\geqslant\check\zeta_{\varepsilon,\varphi}\mleft(y\mright)-\mu$.
If we get $\zeta_\varepsilon\mleft(\widetilde x\widetilde y\mright)\geqslant\check\zeta_{\varepsilon,\varphi}\mleft(x\mright)+\check\zeta_{\varepsilon,\varphi}\mleft(y\mright)+\frac34-2\mu$, we can conclude by applying \eqref{costofdecompositionint}, \eqref{costofdecompositionfrp} and \eqref{costofdecompositiondfrp}.

Decompose $\widetilde x$ and $\widetilde y$ in \eqref{derhamwittdecomposition}.
Applying propositions \ref{xfrpminoration} and \ref{tfandvfhigherdegree}, we are almost done using the pseudovaluation properties of $\zeta_\varepsilon$ \cite[theorem 3.17]{pseudovaluationsonthederhamwittcomplex}, except for the term $\widetilde x|_\mathrm{int}\widetilde y|_\mathrm{frp}$.
For that last one, we need \cite[proposition 3.8]{pseudovaluationsonthederhamwittcomplex} to
conclude.\end{proof}\begin{enonce}{Lemma}\label{assumegooddelta}Let
$x\in W^\dagger\Omega_{\overline R/k}$.
After shrinking $\delta$ if needed, one can assume that for every
$\varepsilon\in\mleft]0,\delta\mright]$:\begin{align*}\check\zeta_{\varepsilon,\varphi}\mleft(x|_\mathrm{int}\mright)&\geqslant-\frac14\text,\\\check\zeta_{\varepsilon,\varphi}\mleft(x|_\mathrm{frp}\mright)&\geqslant\frac12\text,\\\check\zeta_{\varepsilon,\varphi}\mleft(x|_\mathrm{d(frp)}\mright)&\geqslant\frac34\text.\end{align*}\end{enonce}\begin{proof}Let
$\mu>0$.
Let $\widetilde x\in W^\dagger\Omega^{\underline X}_{k\mleft[\underline X\mright]/k}$ be a preimage of $x$ by $\varphi^*$ such that $\zeta_\varepsilon\mleft(\widetilde x\mright)\geqslant\check\zeta_{\varepsilon,\varphi}\mleft(x\mright)-\mu$.

Recall that $n\in\mathbb N$ is the number of indeterminates in $k\mleft[\underline X\mright]$.
When $n=0$, we are working in the ring of Witt vectors $W\mleft(k\mright)$ so that we only have an integral part, and in that case $\zeta_\varepsilon=0$.

When $n>0$, then by \cite[proposition 5.9]{localstructureoftheoverconvergentderhamwittcomplex} we get by induction on $m\in\mathbb N$ that $\zeta_\frac\varepsilon{\mleft(2n\mright)^m}\mleft(\widetilde x\mright)\geqslant\frac{\zeta_\varepsilon\mleft(\widetilde x\mright)}{\mleft(2n\mright)^m}$, from which we get $\check\zeta_{\varepsilon,\varphi}\mleft(x\mright)\geqslant\zeta_\varepsilon\mleft(x\mright)\geqslant-\frac14$ when $\varepsilon$ is small enough.
So we only have to focus on the fractional part.

As $\zeta_\varepsilon$ is a pseudovaluation \cite[theorem 3.17]{pseudovaluationsonthederhamwittcomplex}, and applying \eqref{frpfromdfrp}, we can assume without loss of generality that $\widetilde x\in W^\dagger\Omega^{\underline X-\mathrm{frp},t}_{k\mleft[\underline X\mright]/k}$, and prove that $\check\zeta_{\varepsilon,\varphi}\mleft(\widetilde x\mright)\geqslant\frac34$.

By definition, for every $e\in G\mleft(t\mright)$ there is $s_e\in W^\dagger\Omega^{\mathrm{int},0}_{k\mleft[\underline X\mright]/k}$ such
that:\begin{align*}\sum_{e\in G\mleft(t\mright)}t_F\mleft(s_e\mright)e&=\widetilde x\text,\\\zeta_\varepsilon\mleft(s_e\mright)+\zeta_\varepsilon\mleft(e\mright)&\geqslant-\frac14\text.\end{align*}

The definition of the elements of $G\mleft(t\mright)$ \cite[proposition 5.9]{localstructureoftheoverconvergentderhamwittcomplex} and of $\zeta_\varepsilon$ \cite[definition 3.4]{pseudovaluationsonthederhamwittcomplex} immediately yield that $\zeta_\varepsilon\mleft(e\mright)\geqslant1-\varepsilon$ for every $e\in G\mleft(t\mright)$.
So for every such $e$, we have $\zeta_\varepsilon\mleft(s_e\mright)\geqslant\varepsilon-\frac54$.
We can use the same trick as above to see that, by dividing $\varepsilon$ as much as needed, we have $\zeta_\varepsilon\mleft(s_e\mright)\geqslant\varepsilon-\frac14$ independently on $e$, and we are
done.\end{proof}\section{Finite projective modules}In
all this section, $\mleft(X,\mathcal O_X\mright)$ denotes a ringed space.
We will study finite projective modules, their relationship with $F$-isocrystals, and how one can roll back to the locally free
setting.\begin{enonce}{Definition}A
$\mathcal O_X$-module $M$ is \emph{globally finite projective} if there exists finite sets $J$ and $L$ and an exact sequence of
$\mathcal O_X$-modules:\begin{equation*}\begin{tikzcd}\bigoplus_{j\in J}\mathcal O_X\arrow[r]&\bigoplus_{l\in L}\mathcal O_X\arrow[r,"\psi"]&M\arrow[r]&0\end{tikzcd}\end{equation*}such
that $\psi$ has a section.

A $\mathcal O_X$-module $M$ is \emph{locally finite projective} if there exists a cover $\mleft\{U_i\mright\}_{i\in I}$ of $X$ by opens, where $I$ is any set, such that for each $i\in I$ the $\mathcal O_X|_{U_i}$-module $M|_{U_i}$ is globally finite
projective.\end{enonce}In
other terms, $M$ is a $\mathcal O_X$-module of finite presentation, for which we can globally (respectively locally) chose a finite presentation which splits.

In many settings, the category of $\mathcal O_X$-modules has no projective objects except for the zero module.
What justifies the terminology chosen here is rather the following
fact.\begin{enonce}{Proposition}\label{locallyprojectivesections}Let
$M$ be a globally finite projective module.
Then, the $\mathcal O_X\mleft(X\mright)$-module $M\mleft(X\mright)$ is finite and projective.
Moreover, we have $M\cong\widetilde{M\mleft(X\mright)}$ as
$\mathcal O_X$-modules.\end{enonce}\begin{proof}This
is really just an immediate consequence of the definition.
Since the map $r$ has a section for all $I$, we get that $M\mleft(X\mright)$ is a direct summand of $\bigoplus_{l\in L}\mathcal O_X\mleft(X\mright)$ in the category of $\mathcal O_X\mleft(X\mright)$-modules.

The last part of the statement is given by the description of $\widetilde{M\mleft(X\mright)}$ in
\cite[01BH]{stacksproject}.\end{proof}Take
care, however, not to confuse these sheaves with finite locally free ones.
On general ringed spaces, these notions are not equivalent.
For instance, if $X$ only has a single point and $\mathcal O_X\mleft(X\mright)=\mathbb Z/6\mathbb Z$, then $\widetilde{\mathbb Z/2\mathbb Z}$ yields a counter-example.
Even in more geometric contexts, it is unclear whether these notions are equivalent; see for instance question \ref{classicalquestion} below.

For a commutative ring $R$, and a topological space $X$ endowed with a sheaf $\mathcal F$ of strictly commutative $R$-dgas, we shall denote by $\mathcal F\textsf-\mathsf{LocFProjMIC}$ the full subcategory of $\mathcal F\textsf-\mathsf{MIC}$ whose objects are integrable $\mathcal F$-connections on a locally finite projective $\mathcal F_0$-module.
Similarly, $\mathcal F\textsf-\mathsf{GlobFProjMIC}$ shall denote the full subcategory of $\mathcal F\textsf-\mathsf{MIC}$ whose objects are integrable $\mathcal F$-connections on a globally finite projective $\mathcal F_0$-module.

Notice that the functors defined in proposition \ref{functorfrommorphismintegrableconnections} restrict in the above categories.

A globally finite projective $\mathcal O_X$-module is a direct summand of a finite free $\mathcal O_X$-module.
As we are about to see, this is also the case for $\mathcal F$-connections on a globally finite projective $\mathcal O_X$-module, but not necessarily in the category $\mathcal F\textsf-\mathsf{GlobFProjMIC}$.

In the remainder of this paper, except where otherwise stated for any function $f\colon E\to F$ we shall denote also by $f\colon\operatorname{Mat}\mleft(E\mright)\to\operatorname{Mat}\mleft(F\mright)$ the function it induces on matrices, by applying $f$ to each coefficients.
We shall also do so for functions using unusual notations, such as
$\bullet|_\mathrm{int}$.\begin{enonce}{Lemma}\label{liftingconnections}Let
$\nabla$ be an $\mathcal F$-connection on a globally finite projective $\mathcal O_X$-module $M$.
Let $\psi\colon\bigoplus_{i=1}^r\mathcal F_0\to M$ be a surjective morphism of $\mathcal F_0$-modules with a section $s$, where $r\in\mathbb N^*$.
Fix a $\mathcal F_0\mleft(X\mright)$-basis $\mathcal B$ of $\psi\colon\bigoplus_{i=1}^r\mathcal F_0\mleft(X\mright)$, and let $P$ be the representative matrix of $s\circ\psi$ for $\mathcal B$.

Then, there exists an $\mathcal F$-connection $\widetilde\nabla$ on $\bigoplus_{i=1}^r\mathcal F_0$ such that $\psi$ and $s$ promote to horizontal morphisms between $\nabla$ and $\widetilde\nabla$.
Moreover, if $\nabla$ is integrable and $N$ denotes the representative matrix of $\widetilde\nabla$ for $\mathcal B$, then one can choose $\widetilde\nabla$ so that its curvature matrix for $\mathcal B$
is:\begin{equation*}d\mleft(P\mright)Pd\mleft(P\mright)\text.\end{equation*}\end{enonce}\begin{proof}Put
$\mleft(b_j\mright)_{j\in\mleft\llbracket1,r\mright\rrbracket}\coloneqq\mathcal B$.
Let $A=\mleft(a_{i,j}\mright)_{i,j\in\mleft\llbracket1,r\mright\rrbracket}\in\operatorname{Mat}_r\mleft(\mathcal F_1\mright)$ be the unique matrix such that for every $i\in\mleft\llbracket1,r\mright\rrbracket$ we
have:\begin{equation*}\mleft(s\otimes\operatorname{Id}_{\mathcal F_1}\mright)\circ\nabla\mleft(\psi\mleft(b_i\mright)\mright)=\sum_{j=1}^rb_j\otimes a_{j,i}\text.\end{equation*}

Notice that $P^2=P$ and $A=PA$.
Moreover, $AP+Pd\mleft(P\mright)=A$ by the Leibniz rule.

Now, put $N\coloneqq A-d\mleft(P\mright)P$.
Thus:\begin{equation*}PN=PA-Pd\mleft(P\mright)P=A-Pd\mleft(P\mright)+PPd\mleft(P\mright)=A\text.\end{equation*}

Also:\begin{equation*}NP+d\mleft(P\mright)=AP-d\mleft(P\mright)P+d\mleft(P\mright)=AP+Pd\mleft(P\mright)=A\text.\end{equation*}

If we denote by $\widetilde\nabla$ the $\mathcal F$-connection on $\bigoplus_{i=1}^r\mathcal F_0$ whose representative matrix is $N$ for $\mathcal B$, these two equalities imply that $\psi$ promotes to a morphism of $\mathcal F$-connections $\widetilde\nabla\to\nabla$, and so do its section $s$ by proposition \ref{connectionevaluation}.

For now on, we assume that $\nabla$ is integrable.
If we denote by $C$ the curvature matrix of $\widetilde\nabla$ for $\mathcal B$, this implies that $PC=CP=0$.
Also, $C=N^2+d\mleft(N\mright)$ by proposition \ref{integrabilitycharacterisation}.
This implies that:\begin{equation*}C=\mleft(1-P\mright)\mleft(N^2+d\mleft(N\mright)\mright)\mleft(1-P\mright)\text.\end{equation*}

Now, the above relations give us $\mleft(1-P\mright)N=N-A=-d\mleft(P\mright)P$ and also $N\mleft(1-P\mright)=N-A+d\mleft(P\mright)=Pd\mleft(P\mright)$.
In particular, $\mleft(1-P\mright)N\mleft(1-P\mright)=0$.

Thus:\begin{multline*}\mleft(1-P\mright)d\mleft(N\mright)\mleft(1-P\mright)\\\begin{aligned}&=d\mleft(\mleft(1-P\mright)N\mleft(1-P\mright)\mright)-d\mleft(1-P\mright)N\mleft(1-P\mright)+\mleft(1-P\mright)Nd\mleft(1-P\mright)\\&=0+d\mleft(P\mright)Pd\mleft(P\mright)+d\mleft(P\mright)Pd\mleft(P\mright)\text.\end{aligned}\end{multline*}

Hence,
$C=2d\mleft(P\mright)Pd\mleft(P\mright)+\mleft(1-P\mright)N^2\mleft(1-P\mright)=d\mleft(P\mright)Pd\mleft(P\mright)$.\end{proof}Fix
a smooth commutative $k$-algebra $\overline R$ and a smooth commutative $W\mleft(k\mright)$-algebra $R$ lifting $\overline R$.
On its weak completion $R^\dagger$, fix a morphism of rings $F\colon R^\dagger\to R^\dagger$ lifting the Frobenius and compatible with $F\colon W\mleft(k\mright)\to W\mleft(k\mright)$ as we did previously.
This yields a morphism of weak formal schemes $F\colon\operatorname{Spff}\mleft(R^\dagger\mright)\to\operatorname{Spff}\mleft(R^\dagger\mright)$ by \cite[1.1.2.]{devissagedes}.
It extends, by the universal property of the $p$-adically complete de Rham complex, to a morphism of sheaves of $W\mleft(k\mright)$-dgas $F\colon\mathsf\Omega_{\operatorname{Spff}\mleft(R^\dagger\mright)/W\mleft(k\mright)}\mleft[\frac1p\mright]\to\mathsf\Omega_{\operatorname{Spff}\mleft(R^\dagger\mright)/W\mleft(k\mright)}\mleft[\frac1p\mright]$ on
$\operatorname{Spec}\mleft(\overline R\mright)$.
It also yields a morphism of sheaves of $W\mleft(k\mright)$-dgas on
$\operatorname{Spec}\mleft(\overline R\mright)$:\begin{equation*}t_F\colon\mathsf\Omega_{\operatorname{Spff}\mleft(R^\dagger\mright)/W\mleft(k\mright)}\mleft[\frac1p\mright]\to W^\dagger\Omega_{\operatorname{Spec}\mleft(\overline R\mright)/k}\mleft[\frac1p\mright]\text.\end{equation*}

Let $\phi\colon W^\dagger\Omega_{\operatorname{Spec}\mleft(\overline R\mright)/k}\mleft[\frac1p\mright]\to W^\dagger\Omega_{\operatorname{Spec}\mleft(\overline R\mright)/k}\mleft[\frac1p\mright]$ defined on the gradation by the relation $\phi|_{W^\dagger\Omega^i_{\operatorname{Spec}\mleft(\overline R\mright)/k}\mleft[\frac1p\mright]}\coloneqq p^iF|_{W^\dagger\Omega^i_{\operatorname{Spec}\mleft(\overline R\mright)/k}\mleft[\frac1p\mright]}$ for all $i\in\mathbb N$.
Notice that we have $t_F\circ F=\phi\circ t_F$.
This is a morphism of sheaves of $W\mleft(k\mright)$-dgas by \cite[(1.19)]{derhamwittcohomologyforaproperandsmoothmorphism}.

Of course, all these constructions have a convergent analogue for which we shall use the same notations.
We state the following proposition only for overconvergent de Rham--Witt connections, the case in the classical convergent case being totally
similar.\begin{enonce}{Proposition}\label{liftingisocrystals}Let
$\nabla$ be an integrable $W^\dagger\Omega_{\operatorname{Spec}\mleft(\overline R\mright)/k}\mleft[\frac1p\mright]$-connection on a globally finite projective $W^\dagger\mleft(\mathcal O_{\operatorname{Spec}\mleft(\overline R\mright)}\mright)\mleft[\frac1p\mright]$-module $M$.
Assume that there is an isomorphism $\nabla\cong\phi^*\mleft(\nabla\mright)$ of $W^\dagger\Omega_{\operatorname{Spec}\mleft(\overline R\mright)/k}\mleft[\frac1p\mright]$-connections, and that $W^\dagger\Omega_{\overline R/k}$ has no $p$-torsion.

Then, there exists an integrable $W^\dagger\Omega_{\operatorname{Spec}\mleft(\overline R\mright)/k}\mleft[\frac1p\mright]$-connection on a globally finite projective $W^\dagger\mleft(\mathcal O_{\operatorname{Spec}\mleft(\overline R\mright)}\mright)\mleft[\frac1p\mright]$-module $K$ such that the direct sum $M\oplus K$ is a finite free
$W^\dagger\mleft(\mathcal O_{\operatorname{Spec}\mleft(\overline R\mright)}\mright)\mleft[\frac1p\mright]$-module.\end{enonce}\begin{proof}Let
$r\in\mathbb N^*$, and let $\psi\colon\bigoplus_{i=1}^rW^\dagger\mleft(\overline R\mright)\mleft[\frac1p\mright]\to M\mleft(\operatorname{Spec}\mleft(\overline R\mright)\mright)$ be a surjective morphism of $W^\dagger\mleft(\overline R\mright)\mleft[\frac1p\mright]$-modules.
By proposition \ref{locallyprojectivesections}, $M\mleft(\operatorname{Spec}\mleft(\overline R\mright)\mright)$ is projective, so $\psi$ has a section $s$.

Since we assumed that $\nabla\cong\phi^*\mleft(\nabla\mright)$ as $W^\dagger\Omega_{\operatorname{Spec}\mleft(\overline R\mright)/k}\mleft[\frac1p\mright]$-connections, we also have an isomorphism $\varphi\colon M\mleft(\operatorname{Spec}\mleft(\overline R\mright)\mright)\to{M\mleft(\operatorname{Spec}\mleft(\overline R\mright)\mright)}_F$ of $W^\dagger\mleft(\overline R\mright)\mleft[\frac1p\mright]$-modules, where ${M\mleft(\operatorname{Spec}\mleft(\overline R\mright)\mright)}_F$ denotes the extension of scalars of $M\mleft(\operatorname{Spec}\mleft(\overline R\mright)\mright)$ by the Frobenius morphism $F\colon W^\dagger\mleft(\overline R\mright)\mleft[\frac1p\mright]\to W^\dagger\mleft(\overline R\mright)\mleft[\frac1p\mright]$.
Similarly, we shall also denote by $\psi_F\colon\bigoplus_{i=1}^rW^\dagger\mleft(\overline R\mright)\mleft[\frac1p\mright]\to{M\mleft(\operatorname{Spec}\mleft(\overline R\mright)\mright)}_F$ the extension of scalars of $\psi$ by $F$.

Next, we derive from $\psi$ another surjective morphism of $W^\dagger\mleft(\overline R\mright)\mleft[\frac1p\mright]$-modules $\widetilde\psi\colon\bigoplus_{i=1}^{2r}W^\dagger\mleft(\overline R\mright)\mleft[\frac1p\mright]\to M\mleft(\operatorname{Spec}\mleft(\overline R\mright)\mright)$ by sending every new direct summand to zero.
It has a section $\widetilde s$.
By Schanuel's lemma, we have the following isomorphisms of
$W^\dagger\mleft(\overline R\mright)\mleft[\frac1p\mright]$-modules:\begin{equation*}\ker\mleft(\widetilde\psi\mright)\cong\bigoplus_{i=1}^rW^\dagger\mleft(\overline R\mright)\mleft[\frac1p\mright]\oplus\ker\mleft(\psi\mright)\cong\bigoplus_{i=1}^rW^\dagger\mleft(\overline R\mright)\mleft[\frac1p\mright]\oplus\ker\mleft(\psi_F\mright)\cong\ker\mleft(\widetilde\psi\mright)_F\text.\end{equation*}

Let $K\coloneqq\ker\mleft(\widetilde\psi\mright)$.
Denote by $\iota\colon K\to K_F$ the isomorphism we get by composing the three above ones.
We thus have an isomorphism $\varphi\oplus\iota$ of finite free $W^\dagger\mleft(\overline R\mright)\mleft[\frac1p\mright]$-modules.

Let us consider a $W^\dagger\mleft(\overline R\mright)\mleft[\frac1p\mright]$-basis $\mathcal B=\mleft(b_i\mright)_{i\in\mleft\llbracket1,2r\mright\rrbracket}$ of $M\mleft(\operatorname{Spec}\mleft(\overline R\mright)\mright)\bigoplus K$.
Let $P=\mleft(p_{i,j}\mright)_{i,j\in\mleft\llbracket1,2r\mright\rrbracket}$ be the representative matrix of $\widetilde s\circ\widetilde\psi$ for $\mathcal B$.
By \cite[proposition V 5.1]{algebrecommutativemethodesconstructives}, the representative matrix of ${\widetilde s}_F\circ{\widetilde\psi}_F$ for $\mleft(b_i\otimes1\mright)_{i\in\mleft\llbracket1,2r\mright\rrbracket}$ is $F\mleft(P\mright)$, that is, the matrix whose coefficients are the images through $F$ of the ones of $P$.

Notice that we have the following commutative diagram, in which every horizontal line is an
isomorphism:\begin{equation*}\begin{tikzcd}M\mleft(\operatorname{Spec}\mleft(\overline R\mright)\mright)\bigoplus K\arrow[r,"\varphi\oplus\iota"]\arrow[d,"\widetilde\psi"]&{M\mleft(\operatorname{Spec}\mleft(\overline R\mright)\mright)}_F\bigoplus K_F\arrow[d,"{\widetilde\psi}_F"]\\M\mleft(\operatorname{Spec}\mleft(\overline R\mright)\mright)\arrow[r,"\varphi"]\arrow[d,"\widetilde s"]&{M\mleft(\operatorname{Spec}\mleft(\overline R\mright)\mright)}_F\arrow[d,"{\widetilde s}_F"]\\M\mleft(\operatorname{Spec}\mleft(\overline R\mright)\mright)\bigoplus K\arrow[r,"\varphi\oplus\iota"]&{M\mleft(\operatorname{Spec}\mleft(\overline R\mright)\mright)}_F\bigoplus K_F\text.\end{tikzcd}\end{equation*}

For every $i\in\mleft\llbracket1,2r\mright\rrbracket$ we
have:\begin{equation*}\begin{split}\sum_{j=1}^{2r}p_{j,i}\mleft(b_j\otimes1\mright)&=\mleft(\varphi\oplus\iota\mright)^{-1}\mleft(\sum_{j=1}^{2r}p_{j,i}\mleft(\varphi\oplus\iota\mright)\mleft(b_j\otimes1\mright)\mright)\\&=\mleft(\mleft(\varphi\oplus\iota\mright)^{-1}\circ\widetilde s\circ\widetilde\psi\circ\mleft(\varphi\oplus\iota\mright)\mright)\mleft(b_j\otimes1\mright)\\\sum_{j=1}^{2r}p_{j,i}\mleft(b_j\otimes1\mright)&=\mleft({\widetilde s}_F\circ{\widetilde\psi}_F\mright)\mleft(b_j\otimes1\mright)\text.\end{split}\end{equation*}

This means that $P=F\mleft(P\mright)$.
In particular, for every integer $n\in\mathbb N$ we have $d\mleft(P\mright)=p^nF^n\mleft(d\mleft(P\mright)\mright)$ by \cite[(1.19)]{derhamwittcohomologyforaproperandsmoothmorphism}.
This implies that the coefficients of the matrix $d\mleft(P\mright)$ are in $\bigcap_{m\in\mathbb N}p^mW^\dagger\Omega_{\overline R/k}$ because we assumed that $W^\dagger\Omega_{\overline R/k}$ had no $p$-torsion.
Thus, $d\mleft(P\mright)=0$.

Therefore, by lemma \ref{liftingconnections}, there exists an integrable $W^\dagger\Omega_{\operatorname{Spec}\mleft(\overline R\mright)/k}\mleft[\frac1p\mright]$-connection $\widetilde\nabla$ on $\bigoplus_{i=1}^{2r}W^\dagger\mleft(\mathcal O_{\operatorname{Spec}\mleft(\overline R\mright)}\mright)\mleft[\frac1p\mright]$ such that $\widetilde\psi$ and $\widetilde s$, seen here as morphisms of $W^\dagger\mleft(\mathcal O_{\operatorname{Spec}\mleft(\overline R\mright)}\mright)\mleft[\frac1p\mright]$-modules, promote to horizontal morphisms between $\widetilde\nabla$ and $\nabla$.
To conclude, notice that $\widetilde\nabla$ restricts to an integrable $W^\dagger\Omega_{\operatorname{Spec}\mleft(\overline R\mright)/k}\mleft[\frac1p\mright]$-connection on
$\widetilde K$.\end{proof}\section{The functor}We
shall denote here by $X$ a scheme smooth over $k$, and by $\mathcal O_{X/W\mleft(k\mright)}$ the structure sheaf of the crystalline site $\mathsf{Cris}\mleft(X/W\mleft(k\mright)\mright)$.
In this section, we study Etesse's functor $\acute E\colon\mathsf{Cris}\mleft(X/W\mleft(k\mright)\mright)\to W\Omega_{X/k}\textsf-\mathsf{MIC}$ from crystals over $X/W\mleft(k\mright)$ to integrable $W\Omega_{X/k}$-connections \cite[définition 1.2.5]{complexedederhamwittacoefficientsdansuncristal}.
We shall see here how it yields a functor on convergent isocrystals on $X$.

First, we describe it
locally.\begin{enonce}{Proposition}\label{localetesse}Assume
that $X=\operatorname{Spec}\mleft(\overline R\mright)$ is finite étale over $\mathbb A^n_k$ for some $n\in\mathbb N$.
Let $R$ be a smooth commutative $W\mleft(k\mright)$-algebra lifting $\overline R$.
Then, the restriction of $\acute E$ to crystals in finite $\mathcal O_{X,W\mleft(k\mright)}$-modules factors in the following essentially commutative diagram, in which $\iota$ is fully
faithful:\begin{equation*}\begin{tikzcd}\mathsf{Cris_{fin}}\mleft(X/W\mleft(k\mright)\mright)\arrow[d,"\iota"]\arrow[r,"\acute E"]&W\Omega_{X/k}\textsf-\mathsf{MIC}\\\mathsf\Omega_{\operatorname{Spf}\mleft(\widehat R\mright)/W\mleft(k\mright)}\textsf-\mathsf{MIC}\arrow[ru,"{t_F}^*"]\end{tikzcd}\end{equation*}\end{enonce}\begin{proof}The
functor $\iota$ is given by \cite[proposition 1.3.3]{theoriededieudonnecristallineiii}.
For each $n\in\mathbb N^*$ both $R/p^nR$ and $W_n\mleft(\overline R\mright)$ yield divided power thickenings of $\overline R$.
Hence, given a crystal $A$ over $X/W\mleft(k\mright)$, the morphism $t_F$ induces an isomorphism of $W\mleft(\overline R\mright)$-modules $A\mleft(R/p^nR\mright)\otimes_{R/p^nR}W_n\mleft(\overline R\mright)\cong A\mleft(W_n\mleft(\overline R\mright)\mright)$.

In fact, \cite[corollaire 1.3.2]{theoriededieudonnecristallineiii} and \cite[proposition 1.1.6]{complexedederhamwittacoefficientsdansuncristal} explain how we get an integrable $\mathsf\Omega_{\operatorname{Spec}\mleft(R/p^nR\mright)/W\mleft(k\mright)}$-connection $\nabla_{R,n}$ on the $\mathcal O_{\operatorname{Spec}\mleft(R/p^nR\mright)}$-module $\widetilde{A\mleft(R/p^nR\mright)}$, and an integrable $W_n\Omega_{X/k}$-connection $\nabla_{W,n}$ on the $\mathcal O_{\operatorname{Spec}\mleft(W_n\mleft(\overline R\mright)\mright)}$-module $\widetilde{A\mleft(W_n\mleft(\overline R\mright)\mright)}$.
From the construction of the above connections, the aforementioned isomorphism actually yields an isomorphism ${t_F}^*\mleft(\nabla_{R,n}\mright)\cong\nabla_{W,n}$ of $W_n\Omega_{X/k}$-connections.
Moreover, by definition $\iota\mleft(A\mright)=\lim_n\nabla_{R,n}$ and $\acute E\mleft(A\mright)=\lim_n\nabla_{W,n}$.

Let $K_n$ fit in the following exact sequence of
$W\mleft(\overline R\mright)$-modules:\begin{equation*}\begin{tikzcd}K_n\to A\mleft(R/p^nR\mright)\otimes_{R}W\Omega_{\overline R/k}/p^nW\Omega_{\overline R/k}\to A\mleft(R/p^nR\mright)\otimes_{R}W_n\Omega_{\overline R/k}\to0\text.\end{tikzcd}\end{equation*}

Using \cite[theorem 4.5]{localstructureoftheoverconvergentderhamwittcomplex}, we see that the limit of the $K_n$ must be naught, so taking the limit of these sequences tells us that $\acute E\mleft(A\mright)$ is $p$-adically complete.

Assume now that $A$ is a crystal in finite $\mathcal O_{X,W\mleft(k\mright)}$-modules.
As $\widehat R$ is a coherent ring, we find that $t_F$ is a flat morphism because $W\mleft(\overline R\mright)$ is a product of copies of $\widehat R$ by \cite[theorem 4.5]{localstructureoftheoverconvergentderhamwittcomplex}.
By base change, so are all its reductions modulo $p^n$.
Thus, by \cite[0912]{stacksproject} we find that ${t_F}^*\mleft(\iota\mleft(A\mright)\mright)\cong\acute E\mleft(A\mright)$ as
$W\Omega_{X/k}$-connections.\end{proof}Recall
that given a lax monoidal functor $F\colon\mathsf C\to\mathsf D$, we have a change of enriching category functor $F_*\colon\mathsf C\textsf-\mathsf{Cat}\to\mathsf D\textsf-\mathsf{Cat}$ from the big category of categories enriched over $\mathsf C$ to the big category of categories enriched over $\mathsf D$.
Given a category $\mathsf K$ enriched over $\mathsf C$, this functor is given by defining a category enriched over $\mathsf D$ with the same objects as $\mathsf K$, and defining $\operatorname{Hom}_{F_*\mleft(\mathsf K\mright)}\mleft(A,B\mright)\coloneqq F\mleft(\operatorname{Hom}_\mathsf K\mleft(A,B\mright)\mright)$ for each objects $A$ and $B$ of $\mathsf K$.

The extension of scalars $\bullet\otimes_{\mathbb Z}\mathbb Q\colon\mathbb Z\textsf-\mathsf{Mod}\to\mathbb Q\textsf-\mathsf{Mod}$ is a lax monoidal functor which we shall denote by $\mathfrak i$.
We use this functor and the above formalism to rephrase the notion of an isogeny category, using the fact that preadditive categories can be seen as categories enriched over
$\mathbb Z\textsf-\mathsf{Mod}$.\begin{enonce}{Definition}Let
$\mathsf C$ be a preadditive category.
The \emph{isogeny category of $\mathsf C$} is the category $\mathfrak i_*\mleft(\mathsf C\mright)$ enriched over
$\mathbb Q\textsf-\mathsf{Mod}$.\end{enonce}Functorially
on $\mathsf C$, we have a functor $\mathsf C\to\mathfrak i_*\mleft(\mathsf C\mright)$ given by the tensorisation of $\mathbb Z\to\mathbb Q$.
This functor is an isomorphism of categories when all the hom-objects of $\mathsf C$ are $\mathbb Q$-vector
spaces.\begin{enonce}{Definition}Let
$\mathsf C$ be a category.
Let $F\colon\mathsf C\to\mathsf C$ be an endofunctor.
A \emph{Frobenius structure} on an object $A$ is an isomorphism $A\cong F\mleft(A\mright)$ in $\mathsf C$.

A morphism between Frobenius structures $i_A$ and $i_B$ respectively on objects $A$ and $B$ of $\mathsf C$ is a morphism $\phi\in\operatorname{Hom}_{\mathsf C}\mleft(A,B\mright)$ such that $F\mleft(\phi\mright)\circ i_A=i_B\circ\phi$.

We will denote by $F\textsf-\mathsf C$ the category of Frobenius structures on objects of $\mathsf C$ with these
morphisms.\end{enonce}In
this article, we shall prove various equivalences of categories between such categories.
Let us give a first one to begin
with.\begin{enonce}{Exemple}\label{fisocrystalsdescription}The
absolute Frobenius endomorphism induces an endomorphism of topoi $F\colon\mathsf{Cris_{fin}}\mleft(X/W\mleft(k\mright)\mright)\to\mathsf{Cris_{fin}}\mleft(X/W\mleft(k\mright)\mright)$, whose pullback yields an endofunctor $F^*$ on the category of crystals in finite $\mathcal O_{X/W\mleft(k\mright)}$-modules.

By \cite[theorem 0.7.2]{theconvergenttoposincharacteristic} and \cite[théorème 2.4.2]{cohomologierigideetcohomologierigideasupportsproprespremierepartie}, we know that the category $F\textsf-\mathsf{Isoc}\mleft(X/W\mleft(k\mright)\mright)$ of $F$-isocrystals on $X$ is equivalent to $\mathfrak i_*\mleft(F^*\mright)\textsf-\mathfrak i_*\mleft(\mathsf{Cris_{fin}}\mleft(X/W\mleft(k\mright)\mright)\mright)$, the category of Frobenius structures on objects of the isogeny category of crystals in finite $\mathcal O_{X,W\mleft(k\mright)}$-modules.
Indeed, by virtue of \cite[(2.3.4)]{cohomologierigideetcohomologierigideasupportsproprespremierepartie}, Berthelot's definition coincides with the one made by Ogus in \cite[definition 2.7]{ogus}.
See \cite[theorem 2.2]{notesonisocrystals} for this modern formulation of the equivalence of
categories.\end{enonce}\begin{enonce}{Proposition}\label{overconvergentdescription}Fix
a smooth commutative $k$-algebra $\overline R$, and a smooth commutative $W\mleft(k\mright)$-algebra $R$ lifting $\overline R$.
Recall that we have previously constructed a morphism of sheaves of $W\mleft(k\mright)$-dgas $F\colon\mathsf\Omega_{\operatorname{Spff}\mleft(R^\dagger\mright)/W\mleft(k\mright)}\to\mathsf\Omega_{\operatorname{Spff}\mleft(R^\dagger\mright)/W\mleft(k\mright)}$ which is compatible with the Frobenius
endomorphisms.

Then, the category $F\textsf-\mathsf{Isoc}^\dagger\mleft(\operatorname{Spec}\mleft(\overline R\mright)/W\mleft(k\mright)\mright)$ of overconvergent $F$-isocrystals on $\operatorname{Spec}\mleft(\overline R\mright)$ is equivalent to the category $F^*\textsf-\mathsf\Omega_{\operatorname{Spff}\mleft(R^\dagger\mright)/W\mleft(k\mright)}\mleft[\frac1p\mright]\textsf-\mathsf{LocFProjMIC}$ of Frobenius structures on integrable $\mathsf\Omega_{\operatorname{Spff}\mleft(R^\dagger\mright)/W\mleft(k\mright)}\mleft[\frac1p\mright]$-connections on a locally finite projective $\mathcal O_{\operatorname{Spff}\mleft(R^\dagger\mright)}\mleft[\frac1p\mright]$-module, and such a module is actually globally finite projective.

In other words, the category $F^*\textsf-\mathsf\Omega_{\operatorname{Spff}\mleft(R^\dagger\mright)/W\mleft(k\mright)}\mleft[\frac1p\mright]\textsf-\mathsf{GlobFProjMIC}$ of Frobenius structures on integrable $\mathsf\Omega_{\operatorname{Spff}\mleft(R^\dagger\mright)/W\mleft(k\mright)}\mleft[\frac1p\mright]$-connections on a globally finite projective $\mathcal O_{\operatorname{Spff}\mleft(R^\dagger\mright)}\mleft[\frac1p\mright]$-module is actually equal to
$F^*\textsf-\mathsf\Omega_{\operatorname{Spff}\mleft(R^\dagger\mright)/W\mleft(k\mright)}\mleft[\frac1p\mright]\textsf-\mathsf{LocFProjMIC}$.\end{enonce}\begin{proof}By
\cite[corollaire 2.5.8]{cohomologierigideetcohomologierigideasupportsproprespremierepartie}, the category $F\textsf-\mathsf{Isoc}^\dagger\mleft(\operatorname{Spec}\mleft(\overline R\mright)/W\mleft(k\mright)\mleft[\frac1p\mright]\mright)$ is equivalent to the category of finite projective $R^\dagger\mleft[\frac1p\mright]$-modules with an integrable connection, in the usual algebraic meaning, and a Frobenius structure.

The functor $\widetilde\bullet$ allows one to define an integrable $\mathsf\Omega_{\operatorname{Spff}\mleft(R^\dagger\mright)/W\mleft(k\mright)}\mleft[\frac1p\mright]$-connection with a Frobenius structure from this data.
Indeed, it commutes with direct sums, so that the module of global sections is the finite projective $R^\dagger\mleft[\frac1p\mright]$-module one started with, and same with the tensor product.
So we only need to show that this functor is essentially surjective and fully faithful.

Given a Frobenius structure on an integrable $\mathsf\Omega_{\operatorname{Spff}\mleft(R^\dagger\mright)/W\mleft(k\mright)}\mleft[\frac1p\mright]$-connection on a locally finite projective $\mathcal O_{\operatorname{Spff}\mleft(R^\dagger\mright)}\mleft[\frac1p\mright]$-module, one can cover $\operatorname{Spff}\mleft(R^\dagger\mright)$ with affine opens $\operatorname{Spff}\mleft(R_i\mright)$ for a finite set $i\in I$, on which the associated module has a finite presentation which splits.
Using proposition \ref{locallyprojectivesections}, we find using \cite[corollaire 2.5.6]{cohomologierigideetcohomologierigideasupportsproprespremierepartie} that we actually have a compatible family of overconvergent $F$-isocrystals on $\operatorname{Spec}\mleft(\overline{R_i}\mright)$, which we can glue by \cite[remarques 2.3.3]{cohomologierigideetcohomologierigideasupportsproprespremierepartie} to an essential preimage of our Frobenius structure.
The full faithfulness works
similarly.\end{proof}As
in the overconvergent case, one has a morphism of sheaves of $W\mleft(k\mright)$-dgas
$F\colon\mathsf\Omega_{\operatorname{Spf}\mleft(\widehat R\mright)/W\mleft(k\mright)}\mleft[\frac1p\mright]\to\mathsf\Omega_{\operatorname{Spf}\mleft(\widehat R\mright)/W\mleft(k\mright)}\mleft[\frac1p\mright]$ on $\operatorname{Spec}\mleft(\overline R\mright)$,.\begin{enonce}{Proposition}\label{convergentdescription}Let
us keep the notations given above.
Then, the category $F\textsf-\mathsf{Isoc}\mleft(\operatorname{Spec}\mleft(\overline R\mright)/W\mleft(k\mright)\mright)$ of convergent $F$-isocrystals on $\operatorname{Spec}\mleft(\overline R\mright)$ is equivalent to the category $F^*\textsf-\mathsf\Omega_{\operatorname{Spf}\mleft(\widehat R\mright)/W\mleft(k\mright)}\mleft[\frac1p\mright]\textsf-\mathsf{LocFProjMIC}$ of Frobenius structures on integrable $\mathsf\Omega_{\operatorname{Spf}\mleft(\widehat R\mright)/W\mleft(k\mright)}\mleft[\frac1p\mright]$-connections on a locally finite projective $\mathcal O_{\operatorname{Spf}\mleft(\widehat R\mright)}\mleft[\frac1p\mright]$-module.

Moreover, the category $F^*\textsf-\mathsf\Omega_{\operatorname{Spf}\mleft(\widehat R\mright)/W\mleft(k\mright)}\mleft[\frac1p\mright]\textsf-\mathsf{LocFProjMIC}$ is actually equal to
$F^*\textsf-\mathsf\Omega_{\operatorname{Spf}\mleft(\widehat R\mright)/W\mleft(k\mright)}\mleft[\frac1p\mright]\textsf-\mathsf{GlobFProjMIC}$.\end{enonce}\begin{proof}This
is the same proof as before, only using a result of Etesse \cite[corollaire 1.2.3]{imagesdirectesii} instead of the one by Berthelot, which he closely
follows.\end{proof}Take
care that, in the previous propositions, it is not clear that locally finite projective is equivalent to finite locally free, as it would be for instance in the classical case of a sheaf of quasi-coherent modules on a scheme.
Actually, after discussing with many experts, which the author thanks here, it seems that the following question is still
open.\begin{enonce}{Question}\label{classicalquestion}With
the above notations, is the category of convergent (respectively overconvergent) $F$-isocrystals on $\operatorname{Spec}\mleft(\overline R\mright)$ equivalent to the category of Frobenius structures on integrable $\mathsf\Omega_{\operatorname{Spf}\mleft(\widehat R\mright)/W\mleft(k\mright)}\mleft[\frac1p\mright]$-connections (respectively integrable $\mathsf\Omega_{\operatorname{Spff}\mleft(R^\dagger\mright)/W\mleft(k\mright)}\mleft[\frac1p\mright]$-connections) on a locally free $\mathcal O_{\operatorname{Spf}\mleft(\widehat R\mright)}\mleft[\frac1p\mright]$-module (respectively a locally free
$\mathcal O_{\operatorname{Spff}\mleft(R^\dagger\mright)}\mleft[\frac1p\mright]$-module)?\end{enonce}Would
the answer to this question be positive, the proofs in the next sections could be simplified.
Without knowing it, we will therefore need to work with projective modules, and not just free ones.
But as we have seen, proposition \ref{liftingisocrystals} allows us to reduce to the case of free modules.

Let $E$ be a $F$-isocrystal on a scheme $X$ smooth over $k$.
The previous discussion explains how applying $\acute E$ yields a Frobenius structure on an integrable $W\Omega_{X/k}\mleft[\frac1p\mright]$-connection on a locally finite projective $W\mleft(\mathcal O_X\mright)\mleft[\frac1p\mright]$-module, for the morphism of $W\mleft(k\mright)$-dgas $\phi\colon W\Omega_{X/k}\mleft[\frac1p\mright]\to W\Omega_{X/k}\mleft[\frac1p\mright]$ that we introduced previously.
Let us describe this locally to explain what we
mean.\begin{enonce}{Proposition}\label{etesseisocrystals}Assume
that $X=\operatorname{Spec}\mleft(\overline R\mright)$ is finite étale over $\mathbb A^n_k$ for some $n\in\mathbb N$.
Let $R$ be a smooth commutative $W\mleft(k\mright)$-algebra lifting $\overline R$.
Then, the isogeny functor $\acute E$ from the category of isocrystals on $X$ factors in the following essentially commutative diagram, in which $\iota$ is an equivalence of
categories:\begin{equation*}\begin{tikzcd}F\textsf-\mathsf{Isoc}\mleft(X/W\mleft(k\mright)\mright)\arrow[d,"\iota"]\arrow[r,"\acute E"]&\phi^*\textsf-W\Omega_{X/k}\mleft[\frac1p\mright]\textsf-\mathsf{GlobFProjMIC}\\F^*\textsf-\mathsf\Omega_{\operatorname{Spf}\mleft(\widehat R\mright)/W\mleft(k\mright)}\mleft[\frac1p\mright]\textsf-\mathsf{GlobFProjMIC}\arrow[ru,"{t_F}^*"]&\end{tikzcd}\end{equation*}\end{enonce}\begin{proof}To
construct this diagram, we shall first consider the composition of functors $\mathsf{Cris_{fin}}\mleft(X/W\mleft(k\mright)\mright)\to W\Omega_{X/k}\textsf-\mathsf{MIC}\to W\Omega_{X/k}\mleft[\frac1p\mright]\textsf-\mathsf{MIC}$, where the first arrow is given by proposition \ref{localetesse}, and the second arrow is given by the localization.
The target category is actually enriched over $\mathbb Q\textsf-\mathsf{Mod}$, so this yields a functor $i_*\mleft(\mathsf{Cris_{fin}}\mleft(X/W\mleft(k\mright)\mright)\mright)\to W\Omega_{X/k}\mleft[\frac1p\mright]\textsf-\mathsf{MIC}$ whose source is the isogeny category of $\mathsf{Cris_{fin}}\mleft(X/W\mleft(k\mright)\mright)$.

Now, the same can be done with $\iota$, but by example \ref{fisocrystalsdescription}, we see that this coincides with the functor described in \cite[corollaire 1.2.3]{imagesdirectesii}.
Conclude with proposition
\ref{convergentdescription}.\end{proof}Our
next goal is to show that $\acute E$ is actually an equivalence of categories, and to describe the essential image of overconvergent $F$-isocrystals on
$X$.\section{Working with de Rham--Witt matrices}In
this section, $\overline R$ denotes a commutative $k$-algebra such that $\operatorname{Spec}\mleft(\overline R\mright)$ is finite étale over $\mathbb A^n_k$ for some $n\in\mathbb N$.

Let $r\in\mathbb N^*$.
We extend the definition of the $p$-adic pseudovaluation $\operatorname v_p$ to matrices in the following
manner:\begin{equation*}\operatorname v_p\colon\begin{array}{rl}\operatorname{Mat}_r\mleft(W^\dagger\Omega_{\overline R/k}\mright)\to&\mathbb R\cup\mleft\{+\infty,-\infty\mright\}\\\mleft(x_{i,j}\mright)_{i,j\in\mleft\llbracket1,r\mright\rrbracket}\mapsto&\min\mleft\{\operatorname v_p\mleft(x_{i,j}\mright)\mid i,j\in\mleft\llbracket1,r\mright\rrbracket\mright\}\text.\end{array}\end{equation*}

For simplicity, we only state the overconvergent version of the results here.
However, all results stated here hold for the usual de Rham--Witt complex.
The difficulty of the proofs are in the control of the overconvergence, and depend on the results of deep theorems given in the author's previous articles.
Needless to say, the reader can easily translate all these results to the convergent case, by removing the daggers and replacing weak completion by $p$-adic completion.
In the proofs, the translating reader only needs to ignore all arguments about overconvergence, and see that all structure theorems here used also have a classical de Rham--Witt
flavour.\begin{enonce}{Proposition}\label{frptimesanything}Let
$x\in W^\dagger\Omega^\mathrm{frp}_{\overline R/k}$ and $y\in W^\dagger\Omega_{\overline R/k}$.
Then, there exists an element $z\in W^\dagger\Omega^\mathrm{frp}_{\overline R/k}$ such that $xy\equiv z\pmod p$.

Put differently, $\operatorname v_p\mleft(\mleft(xy\mright)|_\mathrm{int}\mright)>0$ and
$\operatorname v_p\mleft(\mleft(xy\mright)|_\mathrm{d(frp)}\mright)>0$.\end{enonce}\begin{proof}This
result holds when $\overline R=k\mleft[\underline X\mright]$ and one take the canonical Frobenius lift by \cite[lemma 2.7 and proposition 2.8]{pseudovaluationsonthederhamwittcomplex}.
Because for each $t\in\mathbb N$, we have $G\mleft(t\mright)\subset W^\dagger\Omega^\mathrm{frp}_{k\mleft[\underline X\mright]/k}$ from the definition in \cite[9]{localstructureoftheoverconvergentderhamwittcomplex}, the result holds for any Frobenius lift by \cite[(7.8)]{localstructureoftheoverconvergentderhamwittcomplex} and \cite[IX. §1 proposition 3]{algebrecommutativehuit}.
Using \cite[lemma 9.5]{localstructureoftheoverconvergentderhamwittcomplex}, one can conclude for a general $\overline R$.\end{proof}The
next lemma is a constructive version of
\cite[lemma 5.2]{crystalsandderhamwittconnections}.\begin{enonce}{Lemma}\label{basissplusone}Let
$r,s\in\mathbb N^*$ and let $N\in\operatorname{Mat}_r\mleft(W^\dagger\Omega^1_{\overline R/k}\mright)$ be a matrix such that $N^2+d\mleft(N\mright)=0$ and that $\operatorname v_p\mleft(N|_\mathrm{frac}\mright)\geqslant s$.
Let $U\in\operatorname{Mat}_r\mleft(W^\dagger\Omega^{\mathrm{frp},0}_{\overline R/k}\mright)$ be the unique matrix such that $d\mleft(p^sU\mright)=-N|_\mathrm{d(frp)}$.
Then:\begin{equation*}\operatorname v_p\mleft(\mleft(\mleft(1+p^sU\mright)^{-1}N\mleft(1+p^sU\mright)+\mleft(1+p^sU\mright)^{-1}d\mleft(1+p^sU\mright)\mright)|_\mathrm{frac}\mright)\geqslant s+1\text.\end{equation*}\end{enonce}\begin{proof}Let
us first write $N=E+p^sH$, where $E\in\operatorname{Mat}_r\mleft(t_F\mleft(\mathsf\Omega^1_{R^\dagger/W\mleft(k\mright)}\mright)\mright)$ and $H\in\operatorname{Mat}_r\mleft(W^\dagger\Omega^{\mathrm{frac},1}_{\overline R/k}\mright)$.
We
get:\begin{equation*}\begin{split}E^2+d\mleft(E\mright)&=E^2+d\mleft(N\mright)-p^sd\mleft(H\mright)\\&=E^2-\mleft(E+p^sH\mright)^2-p^sd\mleft(H\mright)\\E^2+d\mleft(E\mright)&=-p^s\mleft(EH+HE+d\mleft(H\mright)\mright)-p^{2s}H^2\text.\end{split}\end{equation*}

We thus get that $d\mleft(H\mright)\equiv-\mleft(EH+HE\mright)|_\mathrm{d(frp)}\pmod p$.
Using proposition \ref{frptimesanything}, we find that $EH|_\mathrm{frp}+H|_\mathrm{frp}E$ is congruent to a matrix with coefficients in $W^\dagger\Omega^{\mathrm{frp},2}_{\overline R/k}$ modulo $p$,
hence:\begin{equation}\label{dhmodp}d\mleft(H\mright)\equiv-\mleft(EH|_\mathrm{d(frp)}+H|_\mathrm{d(frp)}E\mright)|_\mathrm{d(frp)}\pmod p\text.\end{equation}

Moreover, $V\coloneqq\sum_{i\in\mathbb N^*}p^{-s}\mleft(-p^sU\mright)^i$ satisfies $\mleft(1+p^sU\mright)^{-1}=1+p^sV$, and its coefficients are in $W^\dagger\mleft(\overline R\mright)$ according to \cite[lemma 6.4]{localstructureoftheoverconvergentderhamwittcomplex} and \cite[theorem 3.17]{pseudovaluationsonthederhamwittcomplex}.
Also:\begin{multline}\label{explicitbasechangeformula}\mleft(1+p^sV\mright)\mleft(E+p^sH\mright)\mleft(1+p^sU\mright)+\mleft(1+p^sV\mright)d\mleft(1+p^sU\mright)\\=E+p^s\mleft(VE+H+EU+d\mleft(U\mright)\mright)+p^{2s}\mleft(VH+VU+HU+Vd\mleft(U\mright)\mright)+p^{3s}VHU\text.\end{multline}

By proposition \ref{frptimesanything}, $VE+EU$ is congruent to a matrix with coefficients in $W^\dagger\Omega^{\mathrm{frp},1}_{\overline R/k}$ modulo $p$.
We also
have:\begin{equation*}d\mleft(VE+EU\mright)=d\mleft(V\mright)E+Vd\mleft(E\mright)+d\mleft(E\mright)U-Ed\mleft(U\mright)\text.\end{equation*}

By construction, $U\equiv-V\pmod p$ and $d\mleft(U\mright)=-H|_\mathrm{d(frp)}$,
thus:\begin{equation*}d\mleft(V\mright)E-Ed\mleft(U\mright)\equiv H|_\mathrm{d(frp)}E+EH|_\mathrm{d(frp)}\pmod p\text.\end{equation*}

The same argument as above tells us that $Vd\mleft(E\mright)+d\mleft(E\mright)U$ is congruent to a matrix with coefficients in $W^\dagger\Omega^{\mathrm{frp},2}_{\overline R/k}$ modulo $p$,
so:\begin{equation*}d\mleft(VE+EU\mright)\equiv\mleft(H|_\mathrm{d(frp)}E+EH|_\mathrm{d(frp)}\mright)|_\mathrm{d(frp)}\pmod p\text.\end{equation*}

Hence, we derive from \eqref{dhmodp} that $d\mleft(VE+H+EU\mright)\equiv0\pmod p$.
Using proposition \ref{frptimesanything} again, we get $VE+H|_\mathrm{frp}+EU\equiv0\pmod p$.
So in the end we find $VE+H+EU+d\mleft(U\mright)\equiv0\pmod p$, and we conclude with
\eqref{explicitbasechangeformula}.\end{proof}Our
assumption on $\overline R$ implies that $W^\dagger\Omega_{\operatorname{Spec}\mleft(\overline R\mright)/k}$ is $p$-torsion free.
It thus injects into $W^\dagger\Omega_{\operatorname{Spec}\mleft(\overline R\mright)/k}\mleft[\frac1p\mright]$, and the $p$-adic pseudovaluation $\operatorname v_p$ extends naturally in the
localisation.\begin{enonce}{Proposition}\label{almostfullfaithfulness}Let
$\mleft(M,\nabla_M\mright)$ and $\mleft(N,\nabla_N\mright)$ be two $W^\dagger\Omega_{\operatorname{Spec}\mleft(\overline R\mright)/k}\mleft[\frac1p\mright]$-connections on free $W^\dagger\mleft(\mathcal O_{\operatorname{Spec}\mleft(\overline R\mright)}\mright)\mleft[\frac1p\mright]$-modules, and let $\varphi\colon\mleft(M,\nabla_M\mright)\to\mleft(N,\nabla_N\mright)$ be a horizontal morphism.
Let $\mleft\{m_i\mright\}_{i\in I}$ and $\mleft\{n_j\mright\}_{j\in J}$ be respective $W^\dagger\mleft(\overline R\mright)\mleft[\frac1p\mright]$-bases of $M\mleft(\operatorname{Spec}\mleft(\overline R\mright)\mright)$ and $N\mleft(\operatorname{Spec}\mleft(\overline R\mright)\mright)$, where $I$ and $J$ denote sets.
Let $E$, $F$ and $G$ be the respective representative matrices of $\nabla_M$, $\nabla_N$ and $\varphi$ in these bases.
Let $t\in\mathbb N$ be such that $\operatorname v_p\mleft(E\mright)\geqslant t$, $\operatorname v_p\mleft(F\mright)\geqslant t$ and $\operatorname v_p\mleft(G\mright)\geqslant-t$.

Then, if for $s\in\mathbb N$ we have $\operatorname v_p\mleft(E|_\mathrm{frac}\mright)\geqslant s+t$ and $\operatorname v_p\mleft(F|_\mathrm{frac}\mright)\geqslant s+t$, we also have
$\operatorname v_p\mleft(G|_\mathrm{frac}\mright)\geqslant s$.\end{enonce}\begin{proof}Write
$\mleft(e_{l,i}\mright)_{l,i\in I}\coloneqq E$, $\mleft(f_{l,j}\mright)_{l,j\in J}\coloneqq F$ and $\mleft(g_{j,i}\mright)^{j\in J}_{i\in I}\coloneqq G$ for the coefficients of the representative matrices.
That is, we
have:\begin{align*}\forall i\in I,\ \nabla_M\mleft(m_i\mright)&=\sum_{l\in I}m_l\otimes e_{l,i}\text,\\\forall i\in J,\ \nabla_N\mleft(n_i\mright)&=\sum_{j\in J}n_j\otimes f_{j,i}\text,\\\forall i\in I,\ \varphi\mleft(m_i\mright)&=\sum_{j\in J}g_{j,i}n_j\text.\end{align*}

By definition of a horizontal morphism, for each $i\in I$, we
find:\begin{equation*}\sum_{l\in I}\varphi\mleft(m_l\mright)\otimes e_{l,i}=\sum_{j\in J}\mleft(g_{j,i}\nabla_N\mleft(n_j\mright)+n_j\otimes d\mleft(g_{j,i}\mright)\mright)\text.\end{equation*}

In other
words:\begin{equation*}\forall i\in I,\ \sum_{j\in J}n_j\otimes\mleft(\sum_{l\in I}g_{j,l}e_{l,i}\mright)=\sum_{j\in J}n_j\otimes\mleft(d\mleft(g_{j,i}\mright)+\sum_{l\in J}g_{l,i}f_{j,l}\mright)\text.\end{equation*}

Therefore:\begin{equation*}\forall i\in I,\ \forall j\in J,\ d\mleft(g_{j,i}\mright)=\sum_{l\in I}g_{j,l}e_{l,i}-\sum_{l\in J}g_{l,i}f_{j,l}\text.\end{equation*}

Let $u\coloneqq\min\mleft\{\operatorname v_p\mleft(\mleft(g_{j,i}\mright)|_\mathrm{frac}\mright)\mid i\in I,\ j\in J\mright\}$.
We want to prove that $u\geqslant s$.

If $u=+\infty$, there is nothing to prove.
Otherwise, let $i\in I$ and $j\in J$ be such that $\operatorname v_p\mleft(\mleft(g_{j,i}\mright)|_\mathrm{frac}\mright)=u$.
We
have:\begin{multline*}d\mleft(g_{j,i}\mright)=\sum_{l\in I}\mleft(g_{j,i}|_\mathrm{int}e_{l,i}|_\mathrm{int}+g_{j,i}|_\mathrm{frp}e_{l,i}|_\mathrm{int}+g_{j,i}e_{l,i}|_\mathrm{frac}\mright)\\-\sum_{l\in J}\mleft(g_{l,i}|_\mathrm{int}f_{j,l}|_\mathrm{int}+g_{l,i}|_\mathrm{frp}f_{j,l}|_\mathrm{int}+g_{l,i}f_{j,l}|_\mathrm{frac}\mright)\text.\end{multline*}

By proposition \ref{frptimesanything}, for every $l\in I$ we have $\operatorname v_p\mleft(\mleft(g_{j,i}|_\mathrm{frp}e_{l,i}|_\mathrm{int}\mright)|_\mathrm{d(frp)}\mright)\geqslant u+t+1$ and $\operatorname v_p\mleft(\mleft(g_{j,i}e_{l,i}|_\mathrm{frac}\mright)|_\mathrm{d(frp)}\mright)\geqslant s$, and similarly for every $l\in J$ we get the inequalities $\operatorname v_p\mleft(\mleft(g_{l,i}|_\mathrm{frp}f_{j,l}|_\mathrm{int}\mright)|_\mathrm{d(frp)}\mright)\geqslant u+t+1$  and $\operatorname v_p\mleft(\mleft(g_{l,i}f_{j,l}|_\mathrm{frac}\mright)|_\mathrm{d(frp)}\mright)\geqslant s$, so we
conclude.\end{proof}\section{The equivalences of categories}As
in the previous section, $\overline R$ denotes a commutative $k$-algebra such that $\operatorname{Spec}\mleft(\overline R\mright)$ is finite étale over $\mathbb A^n_k$ for some $n\in\mathbb N$.
We let $\varphi\colon k\mleft[\underline X,\underline Y\mright]\to\overline R$ be a surjective morphism of $k\mleft[\underline X\mright]$-algebras.

For every $\varepsilon>0$, we extend the definition of $\check\zeta_{\varepsilon,\varphi}$ to
matrices:\begin{equation*}\check\zeta_{\varepsilon,\varphi}\colon\begin{array}{rl}\operatorname{Mat}_r\mleft(W^\dagger\Omega_{\overline R/k}\mright)\to&\mathbb R\cup\mleft\{+\infty,-\infty\mright\}\\\mleft(x_{i,j}\mright)_{i,j\in\mleft\llbracket1,r\mright\rrbracket}\mapsto&\min\mleft\{\check\zeta_{\varepsilon,\varphi}\mleft(x_{i,j}\mright)\mid i,j\in\mleft\llbracket1,r\mright\rrbracket\mright\}\text.\end{array}\end{equation*}

Let $r\in\mathbb N^*$ and $N\in\operatorname{Mat}_r\mleft(W^\dagger\Omega_{\overline R/k}\mright)$.
By lemma \ref{assumegooddelta}, it makes sense to consider this condition for some fixed
$\varepsilon>0$:\begin{equation}\label{overconvergentcondition}\begin{cases}\check\zeta_{\varepsilon,\varphi}\mleft(N|_\mathrm{int}\mright)\geqslant-\frac14\text,\\\check\zeta_{\varepsilon,\varphi}\mleft(N|_\mathrm{frp}\mright)\geqslant\frac12\text,\\\check\zeta_{\varepsilon,\varphi}\mleft(N|_\mathrm{d(frp)}\mright)\geqslant\frac34\text.\end{cases}\end{equation}\begin{enonce}{Lemma}\label{basechangecontrolsoverconvergence}Let
$\varepsilon>0$.
Let $r\in\mathbb N^*$.
Let $U\in\operatorname{Mat}_r\mleft(V\mleft(W^\dagger\mleft(\overline R\mright)\mright)\mright)$ and let $N\in\operatorname{Mat}_r\mleft(W^\dagger\Omega_{\overline R/k}\mright)$ be two matrices such that $\check\zeta_{\varepsilon,\varphi}\mleft(U\mright)\geqslant\frac34$ and that $N$ satisfies condition \eqref{overconvergentcondition}.

Then, the matrices $\mleft(1+U\mright)N$ and $\mleft(1+U\mright)^{-1}N\mleft(1+U\mright)$ also satisfy condition
\eqref{overconvergentcondition}.\end{enonce}\begin{proof}We
have $NU=\mleft(N|_\mathrm{frp}+N|_\mathrm{d(frp)}\mright)\mleft(U|_\mathrm{int}+U|_\mathrm{frp}\mright)+N|_\mathrm{int}\mleft(U|_\mathrm{int}+U|_\mathrm{frp}\mright)$.
Once the first product is expended, all the products satisfy condition \eqref{overconvergentcondition} by applying \eqref{zetabiggerthancheckzeta}, \eqref{costofdecompositionint}, \eqref{costofdecompositionfrp}, \eqref{costofdecompositiondfrp} and the pseudovaluation properties of $\zeta_\varepsilon$ \cite[theorem 3.17]{pseudovaluationsonthederhamwittcomplex}.

Now, $N|_\mathrm{int}U|_\mathrm{int}$ clearly satisfies condition \eqref{overconvergentcondition} by \cite[theorem 3.17]{pseudovaluationsonthederhamwittcomplex}.
For $N|_\mathrm{int}U|_\mathrm{frp}$, we need the stronger proposition \ref{controllingintfrp}.
We have thus shown that $NU$ satisfies \eqref{overconvergentcondition}.
Permuting $N$ and $U$ yield the same result.

Moreover, we have $\mleft(1+U\mright)^{-1}=\sum_{i\in\mathbb N}\mleft(-U\mright)^i$.
In particular, we find using the same arguments as above
that:\begin{equation*}\check\zeta_{\varepsilon,\varphi}\mleft(\mleft(1+U\mright)^{-1}-1\mright)\geqslant\frac34\text.\end{equation*}

This implies that we can use the same reasoning to
conclude.\end{proof}\begin{enonce}{Lemma}\label{specialproductoverconvergence}Let
$\varepsilon>0$.
Let $r\in\mathbb N^*$.
Let $U\in\operatorname{Mat}_r\mleft(V\mleft(W^\dagger\mleft(\overline R\mright)\mright)\mright)$ such that $\check\zeta_{\varepsilon,\varphi}\mleft(U\mright)\geqslant\frac34$.
Then:\begin{equation*}\check\zeta_{\varepsilon,\varphi}\mleft(\mleft(1+U\mright)^{-1}d\mleft(1+U\mright)\mright)\geqslant\frac34\text.\end{equation*}\end{enonce}\begin{proof}Of
course, $d\mleft(1\mright)=0$.
Recall that $\check\zeta_{\varepsilon,\varphi}\mleft(\mleft(1+U\mright)^{-1}-1\mright)\geqslant\frac34$.
So one concludes with \eqref{zetaepsilonandd} and the above
arguments.\end{proof}We
have seen previously that one can always choose a lift of the Frobenius $F\colon\operatorname{Spff}\mleft(R^\dagger\mright)\to\operatorname{Spff}\mleft(R^\dagger\mright)$ from which we derive a morphism of sheaves of $W\mleft(k\mright)$-dgas $t_F\colon\mathsf\Omega_{\operatorname{Spff}\mleft(R^\dagger\mright)/W\mleft(k\mright)}\mleft[\frac1p\mright]\to W^\dagger\Omega_{\operatorname{Spec}\mleft(\overline R\mright)/k}\mleft[\frac1p\mright]$ on $\operatorname{Spec}\mleft(\overline R\mright)$.
Recall that we also have a morphism of sheaves of $W\mleft(k\mright)$-dgas $\phi\colon W^\dagger\Omega_{\operatorname{Spec}\mleft(\overline R\mright)/k}\mleft[\frac1p\mright]\to W^\dagger\Omega_{\operatorname{Spec}\mleft(\overline R\mright)/k}\mleft[\frac1p\mright]$ which we use for Frobenius structures on
$W^\dagger\Omega_{\operatorname{Spec}\mleft(\overline R\mright)/k}\mleft[\frac1p\mright]$-connections.\begin{enonce}{Proposition}\label{localequivalence}The
following functor, given by proposition \ref{functorfrommorphismintegrableconnections}, going from the category of Frobenius structures on integrable $\mathsf\Omega_{\operatorname{Spff}\mleft(R^\dagger\mright)/W\mleft(k\mright)}\mleft[\frac1p\mright]$-connections on a globally finite projective $\mathcal O_{\operatorname{Spff}\mleft(R^\dagger\mright)}\mleft[\frac1p\mright]$-module to the category of Frobenius structures on integrable $W^\dagger\Omega_{\operatorname{Spec}\mleft(\overline R\mright)/k}\mleft[\frac1p\mright]$-connections on a globally finite projective
$W^\dagger\mleft(\mathcal O_{\operatorname{Spec}\mleft(\overline R\mright)}\mright)\mleft[\frac1p\mright]$-module, is an equivalence of
categories:\begin{equation*}{t_F}^*\colon F^*\textsf-\mathsf\Omega_{\operatorname{Spff}\mleft(R^\dagger\mright)/W\mleft(k\mright)}\mleft[\frac1p\mright]\textsf-\mathsf{GlobFProjMIC}\to\phi^*\textsf-W^\dagger\Omega_{\operatorname{Spec}\mleft(\overline R\mright)/k}\mleft[\frac1p\mright]\textsf-\mathsf{GlobFProjMIC}\text.\end{equation*}\end{enonce}\begin{proof}Start
with a Frobenius structure $\nabla$ on an integrable $W^\dagger\Omega_{\operatorname{Spec}\mleft(\overline R\mright)/k}\mleft[\frac1p\mright]$-connection on a globally finite projective $W^\dagger\mleft(\mathcal O_{\operatorname{Spec}\mleft(\overline R\mright)}\mright)\mleft[\frac1p\mright]$-module $M$.

By proposition \ref{liftingisocrystals}, there exist $r\in\mathbb N$, a morphism of $W^\dagger\mleft(\mathcal O_{\operatorname{Spec}\mleft(\overline R\mright)}\mright)\mleft[\frac1p\mright]$-modules $\psi\colon\bigoplus_{i=1}^rW^\dagger\mleft(\mathcal O_{\operatorname{Spec}\mleft(\overline R\mright)}\mright)\mleft[\frac1p\mright]\to M$, which is surjective with a section $s$, as well as an integrable $W^\dagger\Omega_{\operatorname{Spec}\mleft(\overline R\mright)/k}\mleft[\frac1p\mright]$-connection $\widetilde\nabla$ on $\bigoplus_{i=1}^rW^\dagger\mleft(\mathcal O_{\operatorname{Spec}\mleft(\overline R\mright)}\mright)\mleft[\frac1p\mright]$, which promotes $\psi$ and $s$ to horizontal morphisms between $\widetilde\nabla$ and $\nabla$.

Denote by $\mathcal B_0=\mleft(b_i\mright)_{i\in\mleft\llbracket1,r\mright\rrbracket}$ the canonical $W^\dagger\mleft(\overline R\mright)\mleft[\frac1p\mright]$-basis of $\bigoplus_{i=1}^rW^\dagger\mleft(\overline R\mright)\mleft[\frac1p\mright]$.
Let $N_0$ be the representative matrix of $\widetilde\nabla$ for that basis.
Recall that $\phi^*\mleft(\widetilde\nabla\mright)$ is also a $W^\dagger\Omega_{\operatorname{Spec}\mleft(\overline R\mright)/k}\mleft[\frac1p\mright]$-connection on $\bigoplus_{i=1}^rW^\dagger\mleft(\mathcal O_X\mright)\mleft[\frac1p\mright]$, and its representative matrix for $\mathcal B_0$ is $\phi\mleft(N_0\mright)=pF\mleft(N_0\mright)$.
In particular, after applying $\phi^*$ enough times, we end up with a representative matrix with coefficients in $pW^\dagger\Omega_{\overline R/k}$.

Since we have a Frobenius structure, we have $\nabla\cong\phi^*\mleft(\nabla\mright)$ as $W^\dagger\Omega_{\operatorname{Spec}\mleft(\overline R\mright)/k}\mleft[\frac1p\mright]$-connections.
This means that we could have assumed that $N_0$ had coefficients in $pW^\dagger\Omega_{\overline R/k}$ to start with.

By lemma \ref{assumegooddelta}, we can furthermore assume that $N_0$ satisfies condition \eqref{overconvergentcondition} for a fixed $\varepsilon>0$.
We have the same condition on $\mathcal B_0$, seen as a matrix.

We shall now build by induction on $l\in\mathbb N$ families $\mathcal B_l\in{W^\dagger\mleft(\overline R\mright)\mleft[\frac1p\mright]}^r$ as well as matrices $N_l\in\operatorname{Mat}_r\mleft(W^\dagger\Omega^1_{\overline R/k}\mleft[\frac1p\mright]\mright)$ such
that:\begin{align}\mathcal B_l&\text{ is a }W^\dagger\mleft(\overline R\mright)\mleft[\frac1p\mright]\text{-basis of }\bigoplus_{i=1}^rW^\dagger\mleft(\overline R\mright)\mleft[\frac1p\mright]\text,\label{isabasis}\\N_l&\text{ is the representative matrix of }\widetilde\nabla\text{ for }\mathcal B_l\text,\label{representativeconnectionmatrix}\\\operatorname v_p\mleft(\mathcal B_l\mright)&\geqslant0\text,\label{padicvaluationofbl}\\\operatorname v_p\mleft(N_l\mright)&\geqslant0\text,\label{padicvaluationofnl}\\\mathcal B_l&\text{ satisfies }\eqref{overconvergentcondition}\text{ for }\varepsilon\text,\label{bloverconverges}\\N_l&\text{ satisfies }\eqref{overconvergentcondition}\text{ for }\varepsilon\text,\label{nloverconverges}\\\operatorname v_p\mleft(N_l|_\mathrm{frac}\mright)&>l\text,\label{padicvaluationofnlfrac}\\\mathcal B_l&\equiv\mathcal B_{l+1}\pmod{p^{l+1}}\text.\label{easyb}\\N_l&\equiv N_{l+1}\pmod{p^{l+1}}\text.\label{easyn}\end{align}

We have such data for $l=0$.
Assume that we have it for a given $l\in\mathbb N$.
Put $U\coloneqq1-d^{-1}\mleft(N_l|_\mathrm{d(frp)}\mright)$.
Let $\mathcal B_{l+1}\coloneqq U\mathcal B_l$.
As $U$ is overconvergent and invertible, we get condition \eqref{isabasis}.

Let $N_{l+1}\coloneqq U^{-1}N_lU+U^{-1}d\mleft(U\mright)$.
By proposition \ref{connectionbasechange}, $\mathcal B_{l+1}$ and $N_{l+1}$ also satisfy condition \eqref{representativeconnectionmatrix}.

The induction hypothesis \eqref{padicvaluationofbl} and \eqref{padicvaluationofnl} imply conditions \eqref{padicvaluationofbl} and \eqref{padicvaluationofnl}.

By \eqref{frpfromdfrp} and lemma \ref{basechangecontrolsoverconvergence}, we get condition \eqref{bloverconverges}.

Similarly, and using also lemma \ref{specialproductoverconvergence} we get condition \eqref{nloverconverges}.

Proposition \ref{integrabilitycharacterisation} and the induction hypothesis \eqref{padicvaluationofnlfrac} allow us to apply lemma \ref{basissplusone} and get condition \eqref{padicvaluationofnlfrac}.

Induction hypothesis \eqref{padicvaluationofnlfrac} implies conditions \eqref{easyb} and \eqref{easyn}.

Now that we have these families, by $p$-adic overconvergence given by \eqref{zetabiggerthancheckzeta} we can consider the limit $\mathcal B$ of these $W^\dagger\mleft(\overline R\mright)$-bases, which is also a $W^\dagger\mleft(\overline R\mright)$-basis of $\bigoplus_{i=1}^rW^\dagger\mleft(\overline R\mright)$.
Then, the limit $N$ of the $N_l$ overconverges to the representative matrix of $\widetilde\nabla$ for $\mathcal B$.

But most importantly, $N$ has coefficients in $t_F\mleft(\mathsf\Omega_{R^\dagger/W\mleft(k\mright)}\mright)$.
We can thus apply proposition \ref{almostfullfaithfulness} to the representative matrix $P$ of $s\circ\psi$ for $\mathcal B$ and see that it also has coefficients in $t_F\mleft(\mathsf\Omega_{R^\dagger/W\mleft(k\mright)}\mright)$.
This means that $\widetilde\nabla$ is in the essential image of a $\mathsf\Omega_{\operatorname{Spff}\mleft(R^\dagger\mright)/W\mleft(k\mright)}\mleft[\frac1p\mright]$-connection through the functor ${t_F}^*$, and so does $\nabla$.
Moreover, by proposition \ref{almostfullfaithfulness} again, we also get the full
faithfulness.\end{proof}As
usual, the previous proposition has a convergent analogue, with the same proof except that we do not need to take care of overconvergence.
It is thus natural to ask whether the canonical functor from the category of Frobenius structures on integrable $W^\dagger\Omega_{X/k}$-connections on a locally projective $W^\dagger\mleft(\mathcal O_X\mright)$-module to the category of Frobenius structures on integrable $W\Omega_{X/k}$-connections on a locally projective $W\mleft(\mathcal O_X\mright)$-module is fully faithful, where $X$ is a scheme smooth over $k$.
Indeed, an integral version of this result has been demonstrated by Ertl in the locally free setting
\cite{fullfaithfulnessforoverconvergent}.\begin{enonce}{Proposition}\label{fullfaithfullnessalaertl}Consider
the inclusion $\iota\colon W^\dagger\Omega_{X/k}\mleft[\frac1p\mright]\to W\Omega_{X/k}\mleft[\frac1p\mright]$.
The following functor is fully
faithful:\begin{equation*}\iota^*\colon\phi^*\textsf-W^\dagger\Omega_{X/k}\mleft[\frac1p\mright]\textsf-\mathsf{LocFProjMIC}\to\phi^*\textsf-W\Omega_{X/k}\mleft[\frac1p\mright]\textsf-\mathsf{LocFProjMIC}\text.\end{equation*}\end{enonce}\begin{proof}The
question is local on $X$, so we can assume that $X=\operatorname{Spec}\mleft(\overline R\mright)$ satisfies the conditions introduced at the beginning of this section.
By propositions \ref{overconvergentdescription} and \ref{convergentdescription}, we can also consider the categories where the connections are on globally finite projective modules.

Let $i\colon\mathsf\Omega_{\operatorname{Spff}\mleft(R^\dagger\mright)/W\mleft(k\mright)}\mleft[\frac1p\mright]\to\mathsf\Omega_{\operatorname{Spf}\mleft(\widehat R\mright)/W\mleft(k\mright)}\mleft[\frac1p\mright]$ be the inclusion.
We then have the following essentially commutative
diagram:\begin{equation*}\begin{tikzcd}F^*\textsf-\mathsf\Omega_{\operatorname{Spff}\mleft(R^\dagger\mright)/W\mleft(k\mright)}\mleft[\frac1p\mright]\textsf-\mathsf{GlobFProjMIC}\arrow[r,"i^*"]\arrow[d,"{t_F}^*"]&{\widehat F}^*\textsf-\mathsf\Omega_{\operatorname{Spff}\mleft(R^\dagger\mright)/W\mleft(k\mright)}\mleft[\frac1p\mright]\textsf-\mathsf{GlobFProjMIC}\arrow[d,"{t_{\widehat F}}^*"]\\\phi^*\textsf-W^\dagger\Omega_{\operatorname{Spec}\mleft(\overline R\mright)/k}\mleft[\frac1p\mright]\textsf-\mathsf{GlobFProjMIC}\arrow[r,"\iota^*"]&\phi^*\textsf-W\Omega_{\operatorname{Spec}\mleft(\overline R\mright)/k}\mleft[\frac1p\mright]\textsf-\mathsf{GlobFProjMIC}\end{tikzcd}\end{equation*}

By proposition \ref{localequivalence}, which as we mentioned also applies in the convergent case, both vertical arrows are equivalences of categories.
But by \cite[corollaire 1.2.3]{imagesdirectesii}, the functor $i^*$ is fully faithful because it is equivalent to the functor from overconvergent $F$-isocrystals to convergent
$F$-isocrystals.\end{proof}Notice
that the proof here relies on the fact that the result is equivalent to the full faithfulness of the functor from overconvergent $F$-isocrystals to convergent $F$-isocrystals.
Should we get a proof of this result using only the theory of de Rham--Witt connections as in \cite{fullfaithfulnessforoverconvergent}, we would deduce a new proof of that
theorem.\begin{enonce}{Theorem}\label{equivalenceconvergent}Let
$X$ be a scheme smooth on $k$.
The functor $\acute E$ yields and equivalence of categories between the category $F\textsf-\mathsf{Isoc}\mleft(X/W\mleft(k\mright)\mright)$ of convergent $F$-isocrystals on $X$ and the category $\iota^*\colon\phi^*\textsf-W\Omega_{X/k}\mleft[\frac1p\mright]\textsf-\mathsf{LocFProjMIC}$ of Frobenius structures on integrable $W\Omega_{X/k}\mleft[\frac1p\mright]$-connections on a locally finite projective
$W\mleft(\mathcal O_X\mright)\mleft[\frac1p\mright]$-module.\end{enonce}\begin{proof}The
functor $\acute E$ is of local nature, so we can work locally on $X$.
We can assume that $X$ is of the form $\operatorname{Spec}\mleft(\overline R\mright)$ needed to apply proposition \ref{etesseisocrystals} by \cite[theorem 2]{moreetalecoversofaffinespacesinpositivecharacteristic}.
In that case, this follows from the convergent version of proposition \ref{localequivalence}.
We can see that proposition \ref{etesseisocrystals} allow us to glue, so we
conclude.\end{proof}Now,
if $X=\operatorname{Spec}\mleft(\overline R\mright)$, we can extend the diagram in proposition \ref{etesseisocrystals} in the following
way:\begin{equation*}\begin{tikzcd}[column sep=-40]F\textsf-\mathsf{Isoc}\mleft(X/W\mleft(k\mright)\mright)\arrow[dr]\arrow[rr,"\acute E"]&&\phi^*\textsf-W\Omega_{X/k}\mleft[\frac1p\mright]\textsf-\mathsf{GlobFProjMIC}\\&\widehat F^*\textsf-\mathsf\Omega_{\operatorname{Spf}\mleft(\widehat R\mright)/W\mleft(k\mright)}\mleft[\frac1p\mright]\textsf-\mathsf{GlobFProjMIC}\arrow[ru,"{t_{\widehat F}}^*"]\\&F^*\textsf-\mathsf\Omega_{\operatorname{Spff}\mleft(R^\dagger\mright)/W\mleft(k\mright)}\mleft[\frac1p\mright]\textsf-\mathsf{GlobFProjMIC}\arrow[rd,"{t_F}^*"]\arrow[u,"i^*"]&\\F\textsf-\mathsf{Isoc}^\dagger\mleft(X/W\mleft(k\mright)\mright)\arrow[ur]\arrow[rr]\arrow[uuu]&&\phi^*\textsf-W\Omega_{X/k}\mleft[\frac1p\mright]\textsf-\mathsf{GlobFProjMIC}\arrow[uuu,"\iota^*"]\end{tikzcd}\end{equation*}

By proposition \ref{fullfaithfullnessalaertl} and its proof, all the vertical arrows are fully faithful functors.
The upper triangle is essentially commutative by proposition \ref{etesseisocrystals}, and so is the quadrilateral on the left by \cite[corollaire 1.2.3]{imagesdirectesii}, and the one on the right.
This implies that the whole diagram is essentially commutative.

In particular, by glueing we can only assume that $X$ is a scheme smooth over $k$, and get a functor $\acute E^\dagger\colon F\textsf-\mathsf{Isoc}^\dagger\mleft(X/W\mleft(k\mright)\mright)\to\phi^*\textsf-W^\dagger\Omega_{X/k}\mleft[\frac1p\mright]\textsf-\mathsf{LocFProjMIC}$ from the category of overconvergent isocrystals on $X$ to the category of Frobenius structures on integrable $W^\dagger\Omega_{X/k}\mleft[\frac1p\mright]$-connections on a locally finite projective
$W^\dagger\mleft(\mathcal O_X\mright)\mleft[\frac1p\mright]$-module.\begin{enonce}{Theorem}The
functor $\acute E^\dagger$ is an equivalence of
categories.\end{enonce}\begin{proof}Same
proof as theorem \ref{equivalenceconvergent}, but with the overconvergent versions of our previous
results.\end{proof}These
theorem enable us to reformulate question \ref{classicalquestion} in the following
way:\begin{enonce}{Question}Let
$\nabla\to\phi^*\mleft(\nabla\mright)$ be an object in $\phi^*\textsf-W\Omega_{X/k}\mleft[\frac1p\mright]\textsf-\mathsf{LocFProjMIC}$ (respectively $\phi^*\textsf-W^\dagger\Omega_{X/k}\mleft[\frac1p\mright]\textsf-\mathsf{LocFProjMIC}$).
Is the $W\mleft(\mathcal O_X\mright)\mleft[\frac1p\mright]$-module (respectively the $W^\dagger\mleft(\mathcal O_X\mright)\mleft[\frac1p\mright]$-module) on which $\nabla$ is defined locally
free?\end{enonce}\end{document}